\documentclass[a4paper,12pt,oneside]{amsart}
 
\usepackage{pdfsync}

\usepackage{amsmath,amsfonts,amssymb,amsthm,amscd,latexsym,euscript,xypic}

\usepackage[english]{babel}

\usepackage{comment}
\usepackage{color}
\usepackage{ dsfont }

\newcommand{\Cn}{$C^{(n)}$-}
\newcommand{\VP}{\mathrm{VP}}

\newcommand{\supp}{\text{supp}\,}

\newcommand{\crit}{\mathrm{crit}}
\newcommand{\cof}{\mathrm{cof}}
\newcommand{\dom}{\text{dom}\,}

\DeclareMathOperator{\ord}{ORD}

\newcommand{\Coll}{\mathrm{Coll}}
\newtheorem{theo}{Theorem}[section]
\newtheorem{prop}[theo]{Proposition}
\newtheorem*{theo*}{Theorem}

\newtheorem{claim}{Claim}[theo]
\newtheorem{subclaim}{Subclaim}[claim]
\newtheorem{lemma}[theo]{Lemma}
\newtheorem{cor}[theo]{Corollary}

\theoremstyle{definition} 
\newtheorem{defi}[theo]{Definition}
\newtheorem{notation}[theo]{Notation}

\theoremstyle{remark}
\newtheorem{remark}[theo]{Remark}

\renewcommand{\emptyset}{\varnothing}
\def\s{\subseteq}
\def\forces{\Vdash}

\begin{document}

\title[More on the preservation of large cardinals]{More on the preservation of large cardinals under class forcing}

\author{Joan Bagaria}
\address{ICREA (Instituci\'o Catalana de Recerca i Estudis Avan\c{c}ats) and
\newline \indent Departament de Matem\`atiques i Inform\`atica, Universitat de Barcelona. 
Gran Via de les Corts Catalanes, 585,
08007 Barcelona, Catalonia.}
\email{joan.bagaria@icrea.cat}
\thanks{The research work  of both authors was supported by the Spanish
Government under grant MTM2017-86777-P,
and by the Generalitat de Catalunya (Catalan Government) under grant SGR 270-2017. The research of the second author has also been supported by MECD Grant FPU15/00026.
}
\author{Alejandro Poveda}
\address{Departament de Matem\`atiques i Inform\`atica, Universitat de Barcelona. 
Gran Via de les Corts Catalanes, 585,
08007 Barcelona, Catalonia.}
\email{alejandro.poveda@ub.edu}
\address{Current address of the second author: Hebrew University of Jerusalem, Einstein Institute of Mathematics. Jerusalem 91904, Israel.}

\date{\today }
\subjclass[2000]{Primary: 03Exx. Secondary: 03E35, 03E55.}
\keywords{}

\begin{abstract}
We  prove two general results about the preservation of extendible and $C^{(n)}$-extendible cardinals under a wide  class of forcing iterations (Theorems \ref{preservationCn} and \ref{keytheonewfitting}).  As applications we give  new proofs of  the preservation  of Vop\v{e}nka's Principle and $C^{(n)}$-extendible cardinals under Jensen's iteration  for forcing the GCH \cite{JensenForcing}, previously obtained in \cite{Broo} and \cite{Tsa13}, res\-pectively. We prove that $C^{(n)}$-extendible cardinals  are preserved by forcing with standard Easton-support  iterations for any possible $\Delta_2$-definable  behaviour of the power-set function on regular cardinals. We show that one can force proper class-many disagreements between the universe and HOD with respect to the calculation of successors of regular cardinals, while preserving $C^{(n)}$-extendible cardinals. We also show, assuming the GCH,  that the class forcing iteration of Cummings-Foreman-Magidor for forcing  $\diamondsuit_{\kappa^+}^+$ at every $\kappa$ (\cite{CummingsSquare}) preserves $C^{(n)}$-extendible cardinals. We give an optimal result on the consistency of weak square principles and \Cn-extendible cardinals. In the last section we prove another preservation result for \Cn-extendible cardinals under very general (not necessarily definable or weakly homogeneous) class forcing iterations. As  applications we prove the consistency of \Cn-extendible cardinals with $\rm{V}=\rm{HOD}$,  and also with $\mathrm{GA}$ (the Ground Axiom) plus $\mathrm{V}\neq \mathrm{HOD}$, the latter being a strengthening of  a result from \cite{HRW}.

\end{abstract}

\maketitle

\section{Introduction}

The present paper is a contribution to the long-standing program in set theory of studying the robustness of  strong large-cardinal notions under forcing extensions. Specifically, we are interested in the section of the large-cardinal hierarchy ranging between extendible cardinals and  \textit{Vop\v{e}nka's Principle} ($\VP$). 

In a groundbreaking work,  Richard Laver \cite{Lav} proved that supercompactness, one of the most prominent large-cardinal properties,  can be made indestructible under a wide range of forcing notions. Indeed, given a supercompact cardinal $\kappa$, Laver showed that there is a forcing notion (the \emph{Laver preparation}) that preserves the supercompactness  of $\kappa$ and makes it indestructible under further ${<}\kappa$-directed closed forcing. 

Inspired by the work of Laver, several authors subsequently obtained similar  results for other classical large-cardinal notions. For instance,  Gitik and  Shelah  \cite{GitShel} show that a strong cardinal $\kappa$ can be made indestructible under so-called $\kappa^+$-weakly closed forcing satisfying the Prikry condition;  Hamkins  \cite{HamLott} uses the \emph{lottery preparation} forcing to make various types of large cardinals indestructible under appropriate forcing notions. 
More recently,  Brooke-Taylor \cite{Broo} shows that $\VP$ is indestructible under reverse Easton forcing iterations of increasingly directed-closed forcing notions, without the need for any preparatory forcing.  In the present paper we are concerned with the preservation by forcing of $C^{(n)}$-extendible cardinals. This family of large cardinals was introduced in \cite{BagEtAl} (see also \cite{Bag}) as a strengthening of the classical notion of extendibility and was shown to provide natural  milestones on the road from supercompact cardinals up to $\VP$. Extendible cardinals have experienced a renewed interest after Woodin's proof of the HOD-Dichotomy \cite{WOEM1}. Also, $C^{(n)}$-extendible cardinals have had remarkable applications in category theory and algebraic topology (see \cite{BagEtAl}). Thus, the investigation of  the preservation of such cardinals under forcing is a worthwhile project, which may  lead to further  applications.

\smallskip

Recall (see \cite{Bag}) that, for each $n<\omega$, the class $C^{(n)}$ is the $\Pi_n$-definable closed unbounded proper class of all ordinals $\alpha$ that are $\Sigma_n$-correct, i.e., such that $V_\alpha$ is a $\Sigma_n$-elementary substructure of $\rm{V}$. Also, recall that a cardinal $\kappa$ is \emph{$C^{(n)}$-extendible} if for every $\lambda >\kappa$ there is an ordinal  $\mu$ and an elementary embedding $j:V_\lambda \to V_\mu$ such that $\crit(j)=\kappa$, $j(\kappa)>\lambda$, and $j(\kappa)\in C^{(n)}$.\label{DefCnext}

It turns out that $\mathrm{VP}(\Pi_{n+1})$, namely $\VP$ restricted to classes of structures that are $\Pi_{n+1}$-definable,   is equivalent to the existence of a \Cn-extendible cardinal. Hence $\VP$ is equivalent to the existence of a \Cn-extendible cardinal for each $n\geq 1$ (see \cite{Bag} for details). It is in this sense that \Cn-extendible cardinals are canonical representatives of the large-cardinal hierarchy in the region between the first supercompact cardinal and $\VP$. 

\smallskip

In general, the preservation of very large cardinals by forcing  is a delicate issue since it imposes strong forms of agreement between the ground model and the generic forcing extension. For example, suppose $\kappa\in C^{(n)}$ is inaccessible and $\mathbb{P}$ is a ${<}\kappa$-distributive forcing notion. If $\Vdash_{\mathbb{P}}\text{``$\kappa\in \dot{C}^{(n)}$''}$ then $V_\kappa\prec_{\Sigma_n} V^\mathbb{P}$. 
This underlines the fact that the more correct  a large cardinal is,  the harder it is to preserve its correctness under forcing, and therefore the more fragile it becomes. Indeed, one runs  into trouble when seeking a result akin to Laver's indestructibility for supercompact cardinals for stronger large cardinals  such as extendible.
This phenomenon was first pointed out by Tsaprounis in his Ph.D. thesis \cite{TsaPhD} and it was afterwards extensively studied  in \cite{BagHam}, where the following theorem  illustrates the fragility we just described.
\begin{theo}[\cite{BagHam}]\label{BagHamResult}
Suppose that $V_\kappa\prec_{\Sigma_2} V_\lambda$ and $G\subseteq \mathbb{P}$ is a $\rm{V}$-generic filter for nontrivial strategically ${<}\kappa$-closed forcing $\mathbb{P}\in V_\eta$, where $\eta \leq \lambda$ . Then for every $\theta\geq\eta$,
$$V_\kappa=V[G]_\kappa\nprec_{\Sigma_3} V[G]_\theta.$$
\end{theo}
In particular, if $\kappa$ is an extendible cardinal and $\mathbb{P}$ is any non trivial strategically ${<}\kappa$-closed set forcing notion, then forcing with $\mathbb{P}$ destroys the extendibility of $\kappa$. Moreover, the theorem  implies that there is no hope to obtain indestructibility results for $\Sigma_3$-correct large cardinals.
Thus, if one aims for a general theory of preservation of $C^{(n)}$-extendible cardinals one should concentrate on class forcing notions.  

\medskip

The structure of the paper is as follows: In Section 2 we prove that \Cn-extendible cardinals are uniformly characterizable in a Magidor-like way, i.e., similar to Magidor's characterization of supercompact cardinals. This reinforces the fact that $C^{(n)}$-extendible cardinals are a natural model-theoretic strengthening of supercompactness,  first shown in \cite{BagEtAl} (see also \cite{Bag}).  This  characterization of $C^{(n)}$-extendibi\-li\-ty is used in later sections for carrying out preservation arguments under class forcing. The same characterization  has been independently given by W. Boney in \cite{Bo18}, and also in \cite{BaGiSch} for the virtual forms of higher-level analogs of supercompact cardinals (i.e., $n$-remarkable cardinals) and virtual $C^{(n)}$-extendible cardinals. 

Section 3 is devoted to the analysis of  some reflection properties of class forcing iterations that will be useful in later sections for the study of the preservation of $C^{(n)}$-extendible cardinals. Two key notion are that of \textit{adequate iteration} and \textit{$\mathbb{P}$-$\Sigma_k$-reflecting cardinal}. For the same purpose, in Section 4 we introduce the notion of \textit{$\mathbb{P}$-$\Sigma_k$-super\-compact cardinal} and show how  it relates to $C^{(n)}$-extendible cardinals.

In Section 5 we define the notion of suitable iteration and prove a  general  result about the preservation of $C^{(n)}$-extendible cardinals under a wide class of $\rm{ORD}$-length forcing iterations (cf. Theorem \ref{preservationCn}). The prize we pay for considering such general iterations is that we just prove that any $C^{(n+1)}$-extendible cardinal retains its \Cn-extendibility (cf. Corollary \ref{maincor2}).

  Section 6 is focussed on applications of Theorem \ref{preservationCn}. We give a new proof of Brooke-Taylor's theorem on the preservation of $\VP$ \cite{Broo}. The main advantage with respect to the original proof is that our technique allows for a finer control over the amount of Vop\v{e}nka's Principle that is preserved. 

In Section 7 we concentrate on a more concrete class of forcing iterations we call \emph{fitting} and improve the results obtained in Section 5. Our main result here is Theorem \ref{keytheonewfitting}. 
In Section 8 we give several applications of this theorem. First, we show that $C^{(n)}$-extendible cardinals  are preserved by standard Easton class-forcing iterations for any $\Delta_2$-definable possible behaviour of the power-set function on regular cardinals. This extends the main result of \cite{Tsan}. 
Second, with an eye on Woodin's HOD Conjecture, we explore briefly the connections between $C^{(n)}$-extendible cardinals (and thus also $\rm{VP}$) with the principle $\rm{V}=\rm{HOD}$. In particular, we prove  that it is possible to force a complete disagreement, and in many possible forms, between $\rm{V}$ and $\rm{HOD}$ with respect to the calculation of successors of regular cardinals, while $C^{(n)}$-extendible cardinals are preserved. Third, we show that, assuming the GCH, the class forcing iteration of Cummings-Foreman-Magidor that forces $\diamondsuit_{\kappa^+}^+$ at every $\kappa$ (\cite{CummingsSquare}) preserves $C^{(n)}$-extendible cardinals.
Finally, we prove that \Cn-extendible cardinals are consistent with $\square_{\lambda,\cof(\lambda)}$ for a proper class of singular cardinals $\lambda$. This result is optimal in the sense explained in Subsection 8.4.

In Section 9,  we address the question of the preservation of $C^{(n)}$-extendible cardinals under general (non weakly homogeneous, non definable) suitable iterations. 
As  applications of our analysis  we prove the consistency of \Cn-extendible cardinals with $\rm{V}=\rm{HOD}$, and also with $\mathrm{GA}+\mathrm{V}\neq \mathrm{HOD}$. This latter is a strengthening of a result  of Hamkins, Reitz and Woodin \cite{HRW}.

\section{A Magidor-like characterization of \Cn-extendibility}
We shall prove that \Cn-extendible cardinals (defined in  p.~\pageref{DefCnext} of the Introduction) can be charac\-terized in a  way analogous to  the following  characterization of supercompact cardinals due to Magidor.\footnote{Magidor's original characterization does not require that $\lambda$ and $\bar{\lambda}$ are in $C^{(1)}$, or that for every $a\in V_\lambda$ there is $\bar{a}\in V_{\bar{\lambda}}$ with $j(\bar{a})=a$.}
\begin{theo}[\cite{Mag}]\label{MagidorsTheorem}
For a cardinal $\delta$, the following statements are equivalent:
\begin{enumerate}
\item $\delta$ is a supercompact cardinal.
\item For every $\lambda>\delta$ in $C^{(1)}$ and for every $a\in V_\lambda$, there exist ordinals $\bar{\delta}<\bar{\lambda}<\delta$ and there exist some $\bar{a}\in V_{\bar{\lambda}}$  and an elementary embedding $j:V_{\bar{\lambda}}\longrightarrow V_\lambda$ such that:
\begin{itemize}
\item $\crit(j)=\bar{\delta}$ and $j(\bar{\delta})=\delta$.
\item $j(\bar{a})=a$.
\item ${\bar{\lambda}}\in C^{(1)}$.
\end{itemize}
\end{enumerate}
\end{theo}
The existence of a supercompact cardinal is thus characterized by a  form of   reflection for $\Sigma_1$-correct strata of the universe, for it implies that any $\Sigma_1$-truth (i.e., any $\Sigma_1$ sentence, with parameters, true in  $\rm{V}$) is \textit{captured} (up to some change of parameters) by some level below the supercompact cardinal.  The following notion  generalizes this reflection property  to higher levels of complexity.

\begin{defi}[$\Sigma_n$-supercompact cardinal]
Let $n\geq 1$. 
A cardinal $\delta$ is said to be  \emph{$\Sigma_n$-supercompact}  if for every $\lambda>\delta$ in $C^{(n)}$ and $a\in V_\lambda$, there exist $\bar{\delta}<\bar{\lambda}<\delta$ and $\bar{a}\in V_{\bar{\lambda}}$, and there exists an elementary embedding $j:V_{\bar{\lambda}}\longrightarrow V_\lambda$ such that:
\begin{itemize}
\item $\crit(j)=\bar{\delta}$ and $j(\bar{\delta})=\delta$.
\item $j(\bar{a})=a$.
\item ${\bar{\lambda}}\in C^{(n)}$.
\end{itemize}
\end{defi}


\begin{lemma}\label{SigmanSupercompIsCorrect}
For $n\geq 1$, every $\Sigma_n$-supercompact cardinal  belongs to the class $C^{(n+1)}$.
\end{lemma}
\begin{proof}
We prove the lemma by induction on $n\geq 1$. The base case $n=1$ follows by combining Magidor's theorem (cf. Theorem \ref{MagidorsTheorem}) with the well-known fact that supercompact cardinals are $\Sigma_2$-correct \cite[Proposition 22.3]{Kan}. Thus, we shall assume by induction that  the lemma holds for $n$ and prove it for $n+1$.  

Let $\delta$ be a $\Sigma_{n+1}$-supercompact cardinal, and let $\varphi(x,y)$ be a $\Pi_{n+1}$ formula with  $a\in V_\delta$.  Suppose first that $V_\delta$ satisfies the sentence $\exists x\varphi(x,a)$ and let $b\in V_\delta$ be a witness for it. Since $\delta$ is $\Sigma_n$-supercompact, the induction hypothesis guarantees that $\delta\in C^{(n+1)}$, so that $\varphi(b,a)$ is true.

Conversely, suppose that $\exists x\varphi(x,a)$ is true. Let $\lambda<\delta$ be such that $a\in V_\lambda$. Let   $b$  be a witness for this formula and let $\mu\in C^{(n+2)}\setminus \delta^+$ be  such that $b\in V_\mu$. Then, $V_\mu\models \exists x\varphi(x,a)$. By the $\Sigma_{n+1}$-supercompactness of $\delta$ we may find $\bar{\lambda},\bar{\delta}<\bar{\mu}<\delta$, $\bar{\mu}\in C^{(n+2)}$, $\bar{a}\in V_{\overline{\lambda}}$ and $j:V_{\bar{\mu}}\rightarrow V_\mu$ such that $\crit(j)=\bar{\delta}$, $j(\bar{\delta})=\delta$, $j(\bar{\lambda})=\lambda$  and $j(\bar{a})=a$. In particular, $V_{\bar{\mu}}\models \exists x\varphi(x,\bar{a})$. Notice that, by elementarity, $\bar{\lambda}<\bar{\delta}$, hence $a=j(\bar{a})=\bar{a}$, and thus $V_{\bar{\mu}}\models \exists x\varphi(x,a)$. Finally, since $\bar{\mu}<\delta$ and $\bar{\mu}\in C^{(n+2)}$, our induction hypothesis yields $V_{\bar{\mu}}\prec_{\Sigma_{n+1}}V_\delta$. From this latter assertion it is clear that $V_\delta\models \exists x\varphi(x,a)$, as wanted.
\end{proof}

\begin{theo}\label{MagidorLikecharacteri}
For $n\geq 1$, $\delta$ is a  \Cn-extendible cardinal if and only if $\delta$ is $\Sigma_{n+1}$-supercompact.
\end{theo}

\begin{proof}
Suppose that $\delta$ is \Cn-extendible. Fix any $\lambda>\delta$ in $C^{(n+1)}$ and $a\in V_\lambda$. By \Cn-extendibility, let $\mu>\lambda$ in $C^{(n+1)}$, and let $j: V_\mu\longrightarrow V_\theta$ be such that $\crit(j)=\delta$, $j(\delta)>\mu$ and $j(\delta)\in C^{(n)}$, for some ordinal $\theta$. Notice that $j\upharpoonright V_{\lambda}\in V_\theta$.

\begin{claim}
$V_\theta$ satisfies the following sentence:
\begin{eqnarray}
 \exists \bar{\lambda}<j(\delta)\;\exists \bar{\delta}<\bar{\lambda}\;\exists \bar{a}\in V_{\bar{\lambda}}\;\exists j^*: V_{\bar{\lambda}}\longrightarrow V_{j(\lambda)}\nonumber\\(j^*(\bar{a})=j(a)\;\wedge\; j^*(\bar{\delta})=j(\delta)\;\wedge\; V_{\bar{\lambda}}\prec_{\Sigma_{n+1}} V_{j(\lambda)}).\nonumber
\end{eqnarray}
\end{claim}
\begin{proof}[Proof of claim]
It is sufficient to show that $V_\lambda\prec_{\Sigma_{n+1}} V_{j(\lambda)}$, for then the claim follows as witnessed by $\lambda$, $\delta$, $a$, and $j\upharpoonright V_\lambda$. 

 On the one hand, notice that $V_\delta\prec_{\Sigma_{n+1}} V_\mu$, because \Cn-extendible cardinals are $\Sigma_{n+2}-$correct. By elementarity, this implies  $V_{j(\delta)}\prec_{\Sigma_{n+1}} V_{\theta}$. On the other hand, since $j(\delta)>\mu$ and $j(\delta)\in C^{(n)}$, it is true that $V_\mu\prec_{\Sigma_{n+1}} V_{j(\delta)}$ and thus $V_\mu\prec_{\Sigma_{n+1}} V_{\theta}$.
In addition, since $\mu$ and $\lambda$ were both $\Sigma_{n+1}-$correct, it is the case that $V_\lambda\prec_{\Sigma_{n+1}} V_\mu$. Hence, $V_\lambda \prec_{\Sigma_{n+1}}V_\theta$. Also, by elementarity, $V_{j(\lambda)}\prec_{\Sigma_{n+1}} V_{\theta}$. Combining these two facts, we have that $V_\lambda\prec_{\Sigma_{n+1}} V_{j(\lambda)}$.
\end{proof}
By elementarity, $V_\mu$ satisfies the above sentence with the parameters $j(\delta)$ and $j(\lambda)$ replaced by $\delta$ and $\lambda$, respectively. Hence, 
since $\mu \in C^{(n+1)}$, the  sentence is true in the universe. Since $\lambda$ was arbitrarily chosen in $C^{(n+1)}$, this implies  that $\delta$ is a $\Sigma_{n+1}$-supercompact cardinal.

\medskip

For the converse implication, let $\lambda$ be greater than $\delta$ and let us show that there exists an elementary embedding $j: V_\lambda\longrightarrow V_\theta$, for some ordinal $\theta$, such that $\crit(j)=\delta$, $j(\delta)>\lambda$, and $j(\delta)\in C^{(n)}$. Take $\mu>\lambda$ in $C^{(n+1)}$ and let $\bar{\delta},\bar{\lambda}<\bar{\mu}$ and $j:V_{\bar{\mu}}\longrightarrow V_\mu$ be such that $\crit(j)=\bar{\delta}$, $j(\bar{\delta})=\delta$, $j(\bar{\lambda})=\lambda$, and ${\bar{\mu}}\in C^{(n+1)}$. Now notice that the sentence
\begin{equation}\label{formulaCn+}
\exists \alpha\;\exists j^*: V_{\bar{\lambda}}\longrightarrow V_\alpha\,(\crit(j^*)=\bar{\delta}\;\wedge\; j^*(\bar{\delta})>\bar{\lambda}\;\wedge\; {j^*(\bar{\delta})}\in C^{(n)})
\end{equation}
is $\Sigma_{n+1}$-expressible. Moreover, it is true in $\rm{V}$ as witnessed by $\lambda$ and $j$ because $j(\bar{\delta})=\delta>\bar{\lambda}$ and $\delta\in C^{(n)}$ (cf. Lemma \ref{SigmanSupercompIsCorrect}). Thus, since $V_{\bar{\mu}}$ is $\Sigma_{n+1}$-correct and contains $\bar{\delta}$ and $\bar{\lambda}$, it is also true in $V_{\bar{\mu}}$. By elementarity, $V_\mu$ thinks that the sentence
\begin{eqnarray}
\exists \alpha\;\exists j^*: V_{{\lambda}}\longrightarrow V_\alpha\,(\crit(j^*)={\delta}\;\wedge\; j^*({\delta})>{\lambda}\;\wedge\; {j^*(\delta)}\in C^{(n)})\nonumber.
\end{eqnarray} 
is true. Since $\mu\in C^{(n+1)}$, the above displayed sentence is true in $\rm{V}$ and so $\delta$ is $\lambda$-\Cn-extendible. As $\lambda$ was arbitrarily chosen,  $\delta$ is a \Cn-extendible cardinal.
\end{proof}

\begin{remark}\label{RemarkOnC+extendibility}
\rm{
For $\lambda\in C^{(n)}$, a cardinal $\delta$ is called  $\lambda$-$C^{(n)+}$-extendible  if it is $\lambda$-$C^{(n)}$-extendible,  witnessed by some elementary embedding $j\colon V_\lambda\rightarrow V_\theta$ with  $\theta\in C^{(n)}$. Likewise, $\delta$ is called $C^{(n)+}$-extendible if it is $\lambda$-$C^{(n)+}$-extendible for every  $\lambda>\delta$ in $C^{(n)}$ (see \cite[\S4]{Bag}).  A close inspection of the preceding argument  reveals that $\Sigma_{n+1}$-supercompactness is actually equivalent to $C^{(n)+}$-extendibility. This  gives an alternative proof of \cite[Corollary 3.5]{Tsan}.}
\end{remark}

\begin{cor}
A cardinal is extendible if and only if it is $\Sigma_2$-super\-compact.
\end{cor}
A close inspection of the proof of Theorem \ref{MagidorLikecharacteri}  shows that the following holds:
\begin{cor}
\label{equivalence}
For $n\geq 1$, a cardinal $\delta$ is $C^{(n)}$-extendible if and only if for a proper class of  $\lambda$ in $C^{(n+1)}$, for  every $\alpha < \lambda$ there exist $\bar{\delta}, \bar{\alpha}<{\bar{\lambda}}$ and  an elementary embedding $j:V_{\bar{\lambda}}\longrightarrow V_\lambda$ such that:
\begin{itemize}
\item $\crit(j)=\bar{\delta}$ and $j(\bar{\delta})=\delta$.
\item $j(\bar{\alpha})=\alpha$.
\item ${\bar{\lambda}}\in C^{(n+1)}$.
\end{itemize}
\end{cor}



In the light of the above results it is natural to define the class of \emph{$C^{(0)}$-extendible cardinals} as the class of supercompact cardinals.\label{conventionsupercompact} 
Note that since every $C^{(n+1)}$-extendible cardinal is a limit of $C^{(n)}$-extendible cardinals (see  \cite{Bag}), every $\Sigma_{n+1}$-supercompact cardinal  is a limit of $\Sigma_n$-supercompact cardinals.
It will become apparent in the following sections that the notion of $\Sigma_{n+1}$-supercompactness is a useful reformulation of \Cn-extendibility in the context of class forcing.

\section{Some reflection properties for class forcing iterations}

In the sequel we will only work with ORD-length  forcing iterations, since extendible cardinals are generally destroyed by set-size ones (see Theorem \ref{BagHamResult} and the related discussion).  

\medskip

If $\mathbb{P}$ is a forcing iteration, and $G$ is $\mathbb{P}$-generic over $V$, then for every ordinal $\lambda$,  we denote as customary by $G_\lambda$ the $\mathbb{P}_\lambda$-generic filter induced by $G$, i.e., $G_\lambda:=G\cap \mathbb{P}_\lambda$. Also, as usual we denote by $i_{G_\lambda}(\tau)$ the interpretation of the $\mathbb{P}_\lambda$-name $\tau$ by the filter $G_\lambda$.

For the main preservation results given in the following sections we will need to ensure that there are  many  cardinals $\lambda$ that satisfy $V_\lambda[G_\lambda]=V[G]_\lambda$.

\begin{defi}\label{reflecting}
Let $\mathbb{P}$ be a  forcing iteration. A cardinal $\lambda$ is \emph{$\mathbb{P}$-reflecting} if
 $\mathbb{P}$ forces that $V[\dot{G}]_\lambda= V_\lambda[\dot{G}_\lambda]$. 
\end{defi}

\begin{remark}
\label{remarkPReflecting}
Note that since the rank of $i_{G_\lambda}(\tau)$ in $V[G_\lambda]$ is never bigger than the rank of $\tau$ in $\rm{V}$, for any $\tau\in V^{\mathbb{P}_\lambda}$, we clearly have
$$V_\lambda[G_\lambda]\subseteq V[G_\lambda]_\lambda \subseteq V[G]_\lambda.$$
Thus,  $\mathbb{P}$  always forces $``V_\lambda[\dot{G}_\lambda]\subseteq V[\dot{G}]_\lambda$''. Hence a cardinal $\lambda$ is $\mathbb{P}$-reflecting if and only if $\mathbb{P}$ forces $``V[\dot{G}]_\lambda\subseteq V_\lambda[\dot{G}_\lambda]$''. 
\end{remark}

\begin{prop}\label{Preflectingclosed}
$K:=\{\lambda\mid \text{$\lambda$ is $\mathbb{P}$-reflecting}\}$ is closed.
\end{prop}
\begin{proof}
Let $\kappa$ be an accumulation point of $K$. Let $p\in \mathbb{P}$ and $\tau\in V^\mathbb{P}$ be such that  $p\Vdash_\mathbb{P}\tau\in V[\dot{G}]_\kappa$. Since $\kappa$ is a limit cardinal there is $q\leq_\mathbb{P} p$ and $\lambda<\kappa$ such that $q\Vdash_\mathbb{P} \tau\in V[\dot{G}]_\lambda$. Extending $q$ if necessary, we may find  $\theta\in K\cap \kappa$  above $\lambda$ such that $q\Vdash_\mathbb{P}\tau\in  V[\dot{G}]_\lambda \subseteq V_\theta[\dot{G}_\theta]$. Thus, $\kappa\in K$, as wanted.
\end{proof}

The following proposition gives some sufficient  conditions for a cardinal to be $\mathbb{P}$-reflecting.

\begin{prop}
\label{firstprop}
Suppose $\lambda$ is an inaccessible  cardinal and  $\mathbb{P}$ is a forcing   iteration such that: $\mathbb{P}_\lambda$ is a $\lambda$-cc forcing which preserves the inaccessibility of $\lambda$, $\mathbb{P}_\lambda\subseteq V_\lambda$,    and $\Vdash_{\mathbb{P}_\lambda}\text{``$\,\dot{\mathbb{Q}}$ is $\lambda$-distributive''}$, where $\mathbb{P}\cong \mathbb{P}_\lambda\ast\mathbb{Q}$. Then $\lambda$ is $\mathbb{P}$-reflecting.
\end{prop}

\begin{proof}
%
%
On the one hand, by induction on the rank and using the fact that $\mathbb{P}_\lambda$ is $\lambda$-cc and preserves the inaccessibility of $\lambda$,  one can easily show  that $V[G_\lambda]_\lambda \subseteq V_\lambda [G_\lambda]$.

On the other hand, as $\mathbb{P}_\lambda$ preserves the inaccessibility of $\lambda$,  $|V[G_\lambda]_\lambda|=\lambda$. Hence,  since $\Vdash_{\mathbb{P}_\lambda}\text{``$\,\dot{\mathbb{Q}}$ is $\lambda$-distributive''}$, and so $i_{G_\lambda}(\dot{\mathbb{Q}})$ does not add any new subsets of $V[G_\lambda]_\lambda$, we have  $$V[G]_\lambda \subseteq V[G_\lambda]_\lambda.$$  Hence, $V[G]_\lambda \subseteq V_\lambda[G_\lambda]$.
\end{proof}

\medskip

Let  $\mathcal{L}$ denote the language of set theory augmented with an additional unary predicate $\mathbb{P}$. We will denote by $\Sigma_k^{\mathcal{L}}$ (resp. $\Pi_k^{\mathcal{L}}$) the class of $\Sigma_k$ (resp. $\Pi_k$) formulae of $\mathcal{L}$.
 This choice of  language will be useful when working with expressions involving  a  given iteration $\mathbb{P}$. For instance, we shall need to compute the complexity of  the notion
$$\langle V_\alpha,\in,\mathbb{P}\cap V_\alpha\rangle\prec_{\Sigma_k^{\mathcal{L}}}\langle V,\in,\mathbb{P}\rangle  $$
as a property of $\alpha$, when $\mathbb{P}$ is a definable ORD-length forcing  iteration.  

\begin{remark}\label{remarkcorrectcardinals} 
For such  $\mathbb{P}$ and  $\alpha$,  since we are dealing with iterations it would perhaps seem more natural to consider the predicate $\mathbb{P}_\alpha$ instead of $\mathbb{P}\cap V_\alpha$. However,   if   $\mathbb{P}_\alpha \subseteq V_\alpha$ (and so $\alpha$ is a limit ordinal and  the direct limit is taken at stage $\alpha$ of the iteration),  then  we  have $ \mathbb{P}\cap V_\alpha=\mathbb{P}_\alpha$. 
\end{remark}

It is a well-known fact (see, e.g., \cite[\S13]{Jech}) that the truth predicate (in the language of set theory) for $\Sigma_0$ formulae is $\Delta_1$ definable (i.e., both $\Sigma_1$ and $\Pi_1$-definable); and in general, the truth predicate for $\Sigma_k$ (respectively $\Pi_k$) formulae, for $k\geq 1$, is $\Sigma_k$-definable (respectively $\Pi_k$).  However, if $\mathbb{P}$ is a definable predicate, then the complexity (in the language of set theory) of the truth predicate for formulae in the language $\mathcal{L}$ depends naturally  on the complexity of the definition of $\mathbb{P}$.

\begin{prop}
\label{proptruth}\label{complexitytruthpredicate}
If $\mathbb{P}$ is either $\Sigma_m$ or $\Pi_m$-definable, then the truth predicate $\vDash_{\Sigma_0^{\mathcal{L}}}$ for $\Sigma_0^{\mathcal{L}}$ formulae  is $\Delta_{m+1}$  (i.e., both $\Sigma_{m+1}$ and $\Pi_{m+1}$). In general, for $k\geq 1$, the  truth  predicate $\vDash_{\Sigma_k^{ \mathcal{L}}}$ 
for $\Sigma_k^{\mathcal{L}}$  formulae    (resp.    $\vDash_{\Pi_k^{ \mathcal{L}}}$ for $\Pi_{k}^{\mathcal{L}}$ formulae) is $\Sigma_{m +k}$  (resp. $\Pi_{m+k}$). If $\mathbb{P}$ is $\Delta_m$-definable  (with $m\geq 1$),  then $\vDash_{\Sigma_0^{\mathcal{L}}}$  is $\Delta_m$, and   $\vDash_{\Sigma_k^{ \mathcal{L}}}$ (resp. $\vDash_{\Pi_k^{ \mathcal{L}}}$) is $\Sigma_{m +k-1}$ (resp. $\Pi_{m+k-1}$).
\end{prop}
\begin{proof}
Note first that the only atomic formulas in the language $\mathcal{L}$ are of the form $\text{``$x\in y$''}$, $\text{``$x=y$''}$, or $\text{``$x\in\mathbb{P}\,"$}$, where $x$ and $y$ are variable symbols. Hence, if $\mathbb{P}$ is $\Sigma_m$ (resp. $\Pi_m$) definable, then the formula $\text{``$x\in\mathbb{P}\,"$}$ is equivalent to a $\Sigma_m$ (resp. $\Pi_m$)  formula of the language of set theory. 
It follows that  a Boolean combination of atomic formulas of the language $\mathcal{L}$ is equivalent to a Boolean combination of $\Sigma_m$ and $\Pi_m$ formulas in the language of set theory. Hence, the truth predicate  for  Boolean combinations of atomic formulas in the language $\mathcal{L}$ is $\Delta_{m+1}$ definable. The same applies to formulas with bounded quantifiers. 
By induction on $k$ one can now show (as in \cite[\S13]{Jech}) that $\vDash_{\Sigma_k^{ \mathcal{L}}}$ (resp. $\vDash_{\Pi_k^{ \mathcal{L}}}$) is $\Sigma_{m +k}$ (resp. $\Pi_{m+k}$) definable.

\smallskip

If $\mathbb{P}$ is $\Delta_m$-definable, with $m\geq1$, then ``$x\in \mathbb{P}$" is both $\Sigma_m$ and $\Pi_m$-expressible. It easily follows that $\vDash_{\Sigma_0^{\mathcal{L}}}$  is $\Delta_m$-definable, and 
by induction on $k$ one  readily shows that $\vDash_{\Sigma_k^{ \mathcal{L}}}$ (resp. $\vDash_{\Pi_k^{ \mathcal{L}}}$) is $\Sigma_{m +k-1}$ (resp. $\Pi_{m+k-1}$) definable.
\end{proof}

\begin{defi}
For $k\geq 0$  and a  predicate\footnote{We are interested in the case $\mathbb{P}$ is a class-forcing iteration, but the definition as well as  the next    proposition make sense and hold also for any  predicate $\mathbb{P}$.} $\mathbb{P}$,  we shall denote by $C^{(k)}_\mathbb{P}$ the class of all ordinals $\alpha$ such that $$\text{\,$\langle V_\alpha,\in,\mathbb{P}\cap {V_\alpha}\rangle\prec_{\Sigma_k^{\mathcal{L}}}\langle V,\in,\mathbb{P}\rangle$}.$$
\end{defi}

It is easily seen that  $C^{(k)}_{\mathbb{P}}$ is a closed  unbounded proper class, and   $C^{(k+1)}_\mathbb{P}\subseteq C^{(k)}_\mathbb{P}$, for each $k\geq 0$. Notice that $C^{(0)}_{\mathbb{P}}$ is the class of all ordinals. Let us compute next the complexity of $C^{(k)}_{\mathbb{P}}$, in the language of set theory,  when $\mathbb{P}$ is a definable predicate. 

\begin{notation}\label{NotationGamma}
For $m<\omega$, we will denote by  $\Gamma_m$ (resp.~$\mathbf{\Gamma}_m$) the collection of all the formulas in the language of set theory that are either $\Sigma_m$  or $\Pi_m$ (resp.~with parameters). In the sequel, expressions such as ``$X$ is $\Gamma_m$" should be read as ``$X$ is definable by a formula in $\Gamma_m$".
\end{notation}

\begin{prop}\label{complexityclubclass}
The class $C^{(k)}_\mathbb{P}$ is 
$$
\begin{array}{ll}
\Pi_{m+k},&  \text{if }  k\geq1 \text{ and } \mathbb{P}\mbox{ is }\Gamma_m\\  
\Pi_{m+k-1},&  \text{if } k,m\geq1 \text{ and } \mathbb{P}\mbox{ is }\Delta_m.\\ 
\end{array}
$$
\end{prop}
%

\begin{proof}
The class $C^{(0)}_{\mathbb{P}}$ is $\Pi_0$, being the class of all ordinals.  
If $k\geq 1$, then we have that $\alpha\in C^{(k)}_\mathbb{P}$ if and only if 
$\alpha\in C^{(k-1)}_\mathbb{P}$ and 
$$\forall X,Y (X=V_\alpha \wedge Y=\mathbb{P}\cap X \to \forall \bar{a}\in X\,\forall\varphi \in\Sigma^{\mathcal{L}}_k(\,\vDash_{\Sigma_k^\mathcal{L}}\varphi(\bar{a})\to $$ $$ \langle
X,\in,Y\rangle\vDash \varphi(\bar{a})))$$
Now using induction on $k$ and Proposition \ref{complexitytruthpredicate} the  complexity of  the class $C^{(k)}_{\mathbb{P}}$ is  easily seen to be as desired.
\end{proof}

%
%
%

 \medskip

 It is well-known that if $\mathbb{P}$ is a set-forcing notion, then for every condition $p\in \mathbb{P}$, and $\mathbb{P}$-names $\tau_1,\ldots ,\tau_n$, if $\varphi(x_1,\ldots ,x_n)$ is a $\Sigma_k$ (resp. $\Pi_k$) formula of the language of set theory, with $k\geq 1$, then the sentence ``$p\Vdash_{\mathbb{P}} \varphi (\tau_1,\ldots ,\tau_n)$" is also $\Sigma_k$ (resp. $\Pi_k$), in the parameters $\mathbb{P}$ and $\tau_1,\ldots ,\tau_n$. We shall see next that for   definable  
${\rm{ORD}}$-length forcing iterations   that satisfy some mild conditions  the complexity of the forcing relation for $\mathbb{P}$ depends only on the complexity of the definition of $\mathbb{P}$. Thus, if, e.g., $\mathbb{P}$ is $\Gamma_m$-definable, then the forcing relation for $\Sigma_k$ formulae is $\Sigma_{m+k}$-definable (Proposition \ref{coroformulas})..    

\begin{defi}
\label{defadequate}
Let $\mathbb{P}$ be a definable  ${\rm{ORD}}$-length forcing iteration. Then $\mathbb{P}$ is \emph{adequate}  if  there is a  sufficiently rich finite fragment $ZFC^*$ of ZFC which allows to define the forcing relation $\Vdash_{\mathbb{P}}$ for $\Delta_0$-formulae  of the language of set theory and proving the Forcing Theorem for such formulae, and 
 there is a proper class   of ordinals $\kappa$ such that $V_\kappa \models ZFC^*$ and $\mathbb{P}^{V_\kappa} =\mathbb{P}_\kappa$ (hence also $\mathbb{P}_\kappa =\mathbb{P}\cap V_\kappa$).\footnote{See Remark \ref{remarkcorrectcardinals}.}
\end{defi}

Although there are examples of definable class forcing notions $\mathbb{P}$ for which the forcing relation $\Vdash_{\mathbb{P}}$ is not definable (see \cite[\S7]{HKLNS}), if $\mathbb{P}$ satisfies certain conditions, e.g., is either  \emph{pretame} (see \cite{FRI}),  \emph{progressively-closed} (\cite{Rie}), or \emph{suitable} (Definition \ref{suitableiteration}), then the forcing relation $\Vdash_{\mathbb{P}}$ is definable  and the Forcing Theorem holds for $\mathbb{P}$.

\begin{prop}
\label{lemmabounded}
Let $\mathbb{P}$ be  adequate and let $ZFC^*$ be a  finite fragment of ZFC which is  sufficient for defining the forcing relation $\Vdash_{\mathbb{P}}$ for $\Delta_0$-formulae  of the language of set theory and proving the Forcing Theorem for such formulae. Suppose $M$ is a transitive set and  $\kappa\in  M$ is such that $M_\kappa \models ZFC^\ast$ and $\mathbb{P}^{ M_\kappa} =\mathbb{P}_\kappa$. Then, for every $\Delta_0$-formula $\varphi$ of the language of set theory, every $p\in \mathbb{P}$, and every $\mathbb{P}$-name $\tau$ such that $p,\tau \in M_\kappa$,
$$ p\Vdash_{\mathbb{P}}  \varphi(\tau) \quad \mbox{ iff } \quad    M_\kappa   \models``p\Vdash_{\mathbb{P}_\kappa}\varphi(\tau)".$$
\end{prop}

\begin{proof}
Assume $p\Vdash_{\mathbb{P}}  \varphi(\tau)$. Suppose, for a contradiction, that $q\in \mathbb{P}_\kappa$ is such that $q\leq_{\mathbb{P}_\kappa} p$ and $M_\kappa \models ``q\Vdash_{\mathbb{P}_\kappa} \neg \varphi(\tau)"$.  Suppose $G$ is $\mathbb{P}$-generic over $V$, with $q\in G$. Then
$$V[G]\models  \varphi (i_G(\tau)).$$
Note that, since  $\mathbb{P}^{M_\kappa} =\mathbb{P}_\kappa$, a direct limit is taken at $\kappa$ (for otherwise some conditions in $\mathbb{P}_\kappa$ would be proper classes in $M_\kappa$), and therefore being a dense subset of $\mathbb{P}_\kappa$ is absolute between $M$ and $V$. Thus, $G_\kappa:=G\cap \mathbb{P}_\kappa$  is $\mathbb{P}_\kappa$-generic over $M_\kappa$. Also, since $\tau\in M_\kappa$, hence $\tau\in M_\kappa^{\mathbb{P}_\kappa}$,  $i_{G_\kappa}(\tau)\in M_\kappa [G_\kappa]$ and $i_G(\tau)=i_{G_\kappa}(\tau)$. Since $\varphi$ is $\Delta_0$, hence absolute for transitive classes,
$M_\kappa[G_\kappa] \models \varphi (i_{G_\kappa}(\tau)).$
But as $q\in G_\kappa$, this contradicts   $M_\kappa \models ``q\Vdash_{\mathbb{P}_\kappa}\neg \varphi(\tau)"$.

A similar argument proves the converse.
\end{proof}

\begin{cor}
\label{corobounded}
Let $\mathbb{P}$ be adequate.  Then, for every $\Delta_0$-formula $\varphi$, every $p\in \mathbb{P}$, and every $\mathbb{P}$-name $\tau$, the sentence $``p\Vdash_{\mathbb{P}}  \varphi(\tau)"$ is:
\begin{enumerate}
\item $\Delta_{1}^{\mathcal{L}}$ 
\item $\Delta_{m+1}$, if $\mathbb{P}$ is $\Gamma_{m}$-definable (some $m\geq 1$)
\item $\Delta_{m}$, if $\mathbb{P}$ is $\Delta_{m}$-definable (some $m\geq 1$)
\end{enumerate}
with $p$ and  $\tau$ as parameters,
 \end{cor}
 
 \begin{proof}
  The proof of \ref{lemmabounded} shows that $``p\Vdash_{\mathbb{P}}  \varphi(\tau)"$  holds if and only if 
$$\forall M \, \forall \kappa (  M \mbox{ is a transitive set } \wedge \kappa \in M \wedge M_\kappa\models ZFC^\ast   \wedge  \mathbb{P}^{M_\kappa} =\mathbb{P}_\kappa \,  \wedge$$ $$ p,\tau \in M_\kappa \to    M_\kappa  \models ``p\Vdash_{\mathbb{P}_\kappa}\varphi(\tau)")$$
and also if and only if
$$\exists M  \exists \kappa (M \mbox{ is a transitive set } \wedge \kappa \in M \wedge M_\kappa\models ZFC^\ast   \wedge  \mathbb{P}^{M_\kappa} =\mathbb{P}_\kappa \,  \wedge$$ $$ p,\tau \in M_\kappa \wedge   M_\kappa  \models ``p\Vdash_{\mathbb{P}_\kappa}\varphi(\tau)").$$
Now  the   two displayed sentences above are easily seen to be $\Pi_1^{\mathcal{L}}$ and $\Sigma_1^{\mathcal{L}}$, respectively, with $p$ and $\tau$ as parameters. 

Items (2) and (3) now follow easily from (1) and Proposition \ref{complexitytruthpredicate}.
 \end{proof}

 \begin{prop}
 \label{coroformulas}
 Let $\mathbb{P}$ be   adequate. Let $k\geq 1$. Then, for every $\Sigma_k$ (resp. $\Pi_k$) formula $\varphi(t)$ of the language of set theory, every $p\in \mathbb{P}$, and  every  $\mathbb{P}$-name $\tau$,
 the sentence $``p\Vdash_{\mathbb{P}}  \varphi(\tau)"$ is
 \begin{enumerate}
\item  $\Sigma_{k}^{\mathcal{L}}$ (resp. $\Pi_{k}^{\mathcal{L}}$)
\item $\Sigma_{ m+k}$ (resp. $\Pi_{m+k}$), if $\mathbb{P}$ is $\Gamma_{m}$-definable (some $m\geq 1$)
\item $\Sigma_{ m+k-1}$ (resp. $\Pi_{m+k-1}$), if $\mathbb{P}$ is $\Delta_{m}$-definable (some $m\geq 1$).
\end{enumerate}
\end{prop}

\begin{proof}
First, note that  the class $V^{\mathbb{P}}$ of $\mathbb{P}$-names is $\Delta_{1}^{\mathcal{L}}$. 
Now given a $\Sigma_k$ formula $\exists x \forall y \exists z \ldots \psi (x,y,z,\ldots ,t)$ of the language of set theory, and given a condition $p\in \mathbb{P}$ and   $\tau\in V^{\mathbb{P}}$, 
$$p\Vdash_{\mathbb{P}} \exists x \forall y \exists z \ldots \psi (x,y,z,\ldots ,\tau)$$ if and only if
$$\exists \sigma (\sigma \in V^{\mathbb{P}} \wedge \forall \sigma' (\sigma'\in V^{\mathbb{P}} \to \exists \sigma''(\sigma''\in V^{\mathbb{P}} \wedge \ldots p\Vdash_{\mathbb{P}} \psi (\sigma ,\sigma ', \sigma'', \tau))))$$
Since $p\Vdash_{\mathbb{P}} \psi (\sigma ,\sigma ', \sigma'', \tau)$ is $\Delta_1^{\mathcal{L}}$ (Corollary \ref{corobounded}), the last displayed sentence is $\Sigma_{k}^{\mathcal{L}}$, as wanted.

Items (2) and (3) follow easily from (1) and Proposition \ref{complexitytruthpredicate}.
\end{proof}

Let us compute next the complexity of the notion of $\mathbb{P}$-reflecting cardinal, for adequate $\mathbb{P}$.

\begin{prop}\label{ComplexityPreflecting} 
Let   $\mathbb{P}$ be adequate. Then  the sentence 
``$\kappa$ is a $\mathbb{P}$-reflecting cardinal'' is  
\begin{enumerate}
\item $\Pi_{2}^{\mathcal{L}}$,  with $\kappa$ as a parameter.
\item $\Pi_{m+2}$, if $\mathbb{P}$ is  $\Gamma_m$-definable, with $m\geq 1$.
\item $\Pi_{m+1 }$, if $\mathbb{P}$ is $\Delta_m$-definable, with $m\geq 1$.
\end{enumerate}
 \end{prop}

\begin{proof}
 Note that $\mathbb{P}$ forces  $V[\dot{G}]_\kappa \subseteq V_\kappa [\dot{G}_\kappa]$   if and only if the $\mathcal{L}$-sentence 
$$\forall p, \tau \,(\tau \in V^{\mathbb{P}} \wedge p\Vdash_{\mathbb{P}}``rk(\tau)<\kappa"\to$$
$$ \exists \sigma, q\, (\sigma \in V^{\mathbb{P}}\wedge q\leq p\wedge rk(\sigma)<\kappa\wedge q\Vdash_{\mathbb{P}}``\sigma=\tau"))$$
holds in $\langle V,\in,\mathbb{P}\rangle$. 
Thus, since  $V^{\mathbb{P}}$ is a $\Delta_1^{\mathcal{L}}$-definable class, and the expressions  $``p\Vdash_{\mathbb{P}} rk(\tau)<\kappa"$ and $``q\Vdash_{\mathbb{P}}\sigma=\tau"$ are also $\Delta_1^{\mathcal{L}}$ (Corollary \ref{corobounded}), it is easily seen that  the sentence $\text{``$\kappa$ is a $\mathbb{P}$-reflecting cardinal''}$ is $\Pi_{2}^{\mathcal{L}}$, with $\kappa$ as an additional parameter.

Items (2) and (3) follow easily from (1) using Proposition \ref{complexitytruthpredicate}.
\end{proof}

Let us consider, for the sake of conciseness, the following strengthening of the notion of  $\mathbb{P}$-reflecting cardinal (cf. Definition \ref{reflecting}).

\begin{defi}\label{Psigmanreflecting}
If $k\geq 1$ and $\mathbb{P}$ is an ${\rm{ORD}}$-length forcing   iteration, then  a cardinal $\kappa$ is \emph{$\mathbb{P}$-$\Sigma_k$-reflecting} if it is $\mathbb{P}$-reflecting, belongs to $C^{(k)}_{\mathbb{P}}$, and $\mathbb{P}\cap V_\kappa =\mathbb{P}_\kappa$.
\end{defi}

The next Lemma will be crucial in future arguments. 

\begin{lemma}\label{PsigmareflectingCncardinal}
Suppose   $\mathbb{P}$ is adequate. Then for every $k\geq 0$, if $\kappa$ is $\mathbb{P}$-$\Sigma_k$-reflecting and $V_\kappa\models ZFC^*$,   
then $\mathbb{P}$ forces $V[\dot{G}]_\kappa\prec_{\Sigma_k} V[\dot{G}]$. 
\end{lemma}

\begin{proof}
 The claim  is clear for $k=0$. So, assume inductively that $\mathbb{P}$ forces $V[\dot{G}]_\kappa\prec_{\Sigma_{k-1}} V[\dot{G}]$. Notice that, as $\kappa$ is $\mathbb{P}$-reflecting,  any member of $V[\dot{G}]_\kappa$ is given by a  $\mathbb{P}$-name in $V_\kappa$. Let $\varphi({x})$ be a $\Sigma_k$ formula in the language of set theory and let $\tau\in V_\kappa$ be a $\mathbb{P}_\kappa$-name such that $p\Vdash_{\mathbb{P}}\varphi(\tau)$, for some $p\in\mathbb{P}$.   We only need to show that $p\Vdash_{\mathbb{P}} ``V[\dot{G}]_\kappa \models \varphi (\tau)"$.
 
\begin{claim}
The set of conditions $q\in \mathbb{P}_\kappa$ such that  $q\Vdash_{\mathbb{P}}\varphi(\tau)$ is dense below $p\restriction \kappa$.
\end{claim}

\begin{proof}[Proof of claim]
Suppose, aiming for a contradiction, that $p' \leq_{\mathbb{P}_\kappa} p\restriction \kappa$ is  such that     $q\not\Vdash_{\mathbb{P}}\varphi (\tau)$, for all $q\leq_{\mathbb{P}_\kappa} p'$. Since   
$\langle V_\kappa,\in,\mathbb{P}_\kappa\rangle\prec_{\Sigma_k}\langle V,\in,\mathbb{P}\rangle$, and for every $q\in \mathbb{P}_\kappa$ the sentence ``$q \not \Vdash_{\mathbb{P}}  \varphi(\tau)$" is $\Pi_{k}^{\mathcal{L}}$ (Proposition \ref{coroformulas}), we have that $``q\not \Vdash_{\mathbb{P}_\kappa}  \varphi(\tau)"$ holds in $\langle V_\kappa ,\in , \mathbb{P}_\kappa\rangle$. Therefore, 
$$\langle V_\kappa ,\in , \mathbb{P}_\kappa\rangle\models ``p'  \Vdash_{\mathbb{P}_\kappa} \neg \varphi(\tau)".$$
Again, since  $\langle V_\kappa,\in,\mathbb{P}_\kappa \rangle\prec_{\Sigma_{k}}\langle V,\in,\mathbb{P}\rangle$, and the quoted displayed sentence is $\Pi_{k}^{\mathcal{L}}$,
$$\langle V,\in, \mathbb{P}\rangle \models ``p'   \Vdash_{\mathbb{P}} \neg \varphi(\tau)"$$ which yields  the desired contradiction to the fact that $p\Vdash_{\mathbb{P}}\varphi(\tau)$.
\end{proof}

By the claim,  $p\restriction \kappa \Vdash_{\mathbb{P}}\varphi (\tau)$, and since   $\langle V_\kappa,\in,\mathbb{P}_\kappa\rangle\prec_{\Sigma_{k}}\langle V,\in,\mathbb{P}\rangle$ and the above sentence is $\Sigma^\mathcal{L}_{k}$ we have  $$\langle V_\kappa,\in,\mathbb{P}_\kappa\rangle\vDash ``p\restriction \kappa \Vdash_{\mathbb{P}_\kappa}\varphi (\tau)".$$ 
Thus, if $G$ is $\mathbb{P}$-generic over $V$, with $p\in G$, then $V_\kappa [G_\kappa]\models \varphi (i_{G_\kappa}(\tau))$. Hence, since $i_{G_\kappa}(\tau)=i_G(\tau)$, and $\kappa$ is $\mathbb{P}$-reflecting,    $V[G]_\kappa\models \varphi(i_G(\tau))$, as wanted.
\end{proof}

%

\medskip

For future use, let us calculate the complexity of the predicate ``$\kappa$ is a $\mathbb{P}$-$\Sigma_k$-reflecting cardinal'' for  definable forcing iterations.

\begin{prop}\label{ComplexityPsigmareflecting}
Let  $\mathbb{P}$ be adequate. Then the predicate 
``$\kappa$ is $\mathbb{P}$-$\Sigma_k$-reflecting" is
$$
\begin{array}{ll}
\Pi_{m+2}, & \mbox{if } k=1 \mbox{ and $\mathbb{P}$ is } \Gamma_m\\ 
\Pi_{m+k},  & \mbox{if } k>1 \mbox{ and $\mathbb{P}$ is } \Gamma_m\\
\Pi_{m+1}, & \mbox{if } k=1 \mbox{ and $\mathbb{P}$ is } \Delta_m\\ 
\Pi_{m+k-1},  & \mbox{if } k>1 \mbox{ and $\mathbb{P}$ is } \Delta_m
\end{array}
$$
 \end{prop}

\begin{proof}
First, the assertion $``\mathbb{P}\cap V_\kappa =\mathbb{P}_\kappa"$ is $\Pi_{m+1}$ (in the parameter $\kappa$) if $\mathbb{P}$ is $\Pi_m$-definable, and $\Pi_m$ if $\mathbb{P}$ is $\Sigma_m$-definable. For the quoted sentence holds if and only if:
 $$\forall x (x\in \mathbb{P} \to ({\rm{rk}}(x)<\kappa  \leftrightarrow {\rm{lh}}(x)<\kappa)).$$
Second, by Proposition \ref{ComplexityPreflecting}, the sentence  $\text{``$\kappa$ is a $\mathbb{P}$-reflecting cardinal''}$ is $\Pi_{m+2}$  if $\mathbb{P}$ is $\Gamma_m$-definable; and it is $\Pi_{m+1}$ if $\mathbb{P}$ is $\Delta_m$-definable.
Finally, as shown in Proposition \ref{complexityclubclass},  for $k\geq 0$ the  class $C^{(k)}_\mathbb{P}$   
is $\Pi_{m+k}$  if $\mathbb{P}$ is $\Pi_m$ or $\Sigma_m$-definable; and it is  $\Pi_{m+k-1}$  if $\mathbb{P}$ is $\Delta_m$-definable. 

An easy calculation now shows that the predicate ``$\kappa$ is a $\mathbb{P}$-$\Sigma_k$-reflecting cardinal'' has the claimed  complexity.
\end{proof}
\section{$\mathbb{P}$-$\Sigma_n$-supercompactness}

The following definition gives a refinement  of the notion of $\Sigma_n$-supercom\-pact cardinal, relative to definable iterations.

\begin{defi}[$\mathbb{P}$-$\Sigma_n$-supercompactness]\label{PSigmaSupercom}
If $n\geq 1$ and $\mathbb{P}$ is an ${\rm{ORD}}$-length  forcing  iteration, then   a cardinal $\delta$ is  \emph{$\mathbb{P}$-$\Sigma_n$-supercompact} if there exists a proper class of $\mathbb{P}$-$\Sigma_n$-reflecting cardinals, and for every such cardinal $\lambda>\delta$  and every $a\in V_\lambda$ there exist $\bar{\delta}<\bar{\lambda}<\delta$ and $\bar{a}\in V_{\bar{\lambda}}$, together with an elementary embedding $j: V_{\bar{\lambda}}\longrightarrow V_\lambda$ such that:
\begin{itemize}
\item[(i)] $\crit(j)=\bar{\delta}$ and $j(\bar{\delta})=\delta$
\item[(ii)] $j(\bar{a})=a$
\item[(iii)] $\bar{\lambda}$ is $\mathbb{P}$-$\Sigma_n$-reflecting. 
\end{itemize}
\end{defi}

The next proposition and corollary unveil the connection between the notions of $\Sigma_n$-supercompact and $\mathbb{P}$-$\Sigma_n$-super\-com\-pact cardinals. For  conciseness, let us denote by $\mathcal{S}^{\Sigma_n}$, $\mathcal{S}^{\Sigma_n}_{\mathbb{P}}$, and $\mathcal{E}^{(n)}$ the classes of $\Sigma_n$-supercompact, $\mathbb{P}$-$\Sigma_n$-supercompact, and $C^{(n)}$-extendible cardinals, respectively. 

\begin{prop}\label{PropertiesPsigmasupercom}
Let  $\mathbb{P}$ be an ${\rm{ORD}}$-length forcing iteration, and suppose there is a proper class of $\mathbb{P}$-$\Sigma_n$-reflecting cardinals. Then,

\begin{enumerate}
\item  $\mathcal{S}^{\Sigma_{n}}_{\mathbb{P}}\subseteq \mathcal{S}^{\Sigma_{n}}$

\smallskip
 
\item[(2)] If $\mathbb{P}$ is  adequate and $\Delta_2$-definable, then  $\mathcal{S}^{\Sigma_{3}}\subseteq\mathcal{S}^{\Sigma_1}_\mathbb{P}$.
\smallskip

\item[(3)] If $\mathbb{P}$ is  adequate and $\Delta_{m+1}$-definable, some $m\geq 1$,  and   either  $m$ or $n$ are greater than $1$, then 
 $\mathcal{S}^{\Sigma_{m+n}}\subseteq\mathcal{S}^{\Sigma_n}_\mathbb{P}$. 
 \smallskip

%
\end{enumerate}
In particular, if $\mathbb{P}$ is adequate and    $\Delta_2$, then    every $\Sigma_2$-supercompact cardinal  is $\mathbb{P}$-$\Sigma_1$-supercompact; and for every $n>1$, a cardinal is  $\Sigma_{n}$-supercompact if and only if it is $\mathbb{P}$-$\Sigma_{n}$-supercompact.
\end{prop} 

\begin{proof}
(1): Assume $\delta$ is $\mathbb{P}$-$\Sigma_{n}$-supercompact. Let $\lambda>\delta$ be a $\Sigma_{n}$-correct cardinal and let $a\in V_\lambda$. Let $\kappa>\lambda$ be a $\mathbb{P}$-$\Sigma_n$-reflecting cardinal. Notice that $V_\kappa\vDash``\lambda\in C^{(n)}"$. By $\mathbb{P}$-$\Sigma_n$-supercompactness, there are $\bar{\delta} < \bar{\kappa}<\delta$, with $\bar{\kappa}$ being $\mathbb{P}$-$\Sigma_n$-reflecting,  and there are $\bar{a}\in V_{\bar{\kappa}}$ and $\bar{\lambda} <\bar{\kappa}$, together with an elementary embedding $j:V_{\bar{\kappa}}\longrightarrow V_\kappa$ such that $\crit(j)=\bar{\delta}$  $j(\bar{\delta})=\delta$,  $j(\bar{a})=a$, and $j(\bar{\lambda})=\lambda$. By elementarity, $V_{\bar{\kappa}}$ thinks that $\bar{\lambda}$ is a $\Sigma_{n}$-correct cardinal. Since $\bar{\kappa}\in C^{(n)}$, it follows that $\bar{\lambda}\in C^{(n)}$.

\smallskip

(2) and (3): 
Assume $\delta$ is $\Sigma_3$-supercompact (in case (2), i.e., $m=n=1$) or $\Sigma_{m+n}$-supercompact, in case (3).  Let $\lambda>\delta$ be   $\mathbb{P}$-$\Sigma_{n}$-reflecting, and  let $\kappa>\lambda$ be a   $\Sigma_{m+n}$-correct cardinal. Since  
 being a $\mathbb{P}$-$\Sigma_{n}$-reflecting cardinal  is a $\Pi_{3}$ property (case (2)), and a $\Pi_{m+n}$  property otherwise (case (3)),     $V_\kappa$ thinks that $\lambda$ is  $\mathbb{P}$-$\Sigma_{n}$-reflecting (cf. Proposition \ref{ComplexityPsigmareflecting}). By our assumption, there exist $\bar{\delta}<\bar{\lambda}<\bar{\kappa}<\delta$  with $\bar{\kappa}\in C^{(3)}$ (case (2)), or $\bar{\kappa}\in C^{(m+n)}$, otherwise,  and there exists  an elementary embedding $j:V_{\bar{\kappa}}\longrightarrow V_\kappa$ such that $\crit(j)=\bar{\delta}$,   $j(\bar{\delta})=\delta$,    and $j(\bar{\lambda})=\lambda$. By elementarity, $V_{\bar{\kappa}}$ thinks that $\bar{\lambda}$ is $\mathbb{P}$-$\Sigma_{n}$-reflecting. Since  $V_{\bar{\kappa}}\prec_{\Sigma_{3}}V$ (case (2)), or $V_{\bar{\kappa}}\prec_{\Sigma_{m+n}}V$, otherwise,  $\bar{\lambda}$ is $\mathbb{P}$-$\Sigma_{n}$-reflecting in $\rm{V}$.  Thus, the restriction $j\upharpoonright V_{\bar{\lambda}}$ witnesses the $\mathbb{P}$-$\Sigma_n$-supercompactness of $\delta$. 
\end{proof}

The proposition above together with Theorem \ref{MagidorLikecharacteri} yield the following:

\begin{cor}
\label{coroscext}
Suppose  $\mathbb{P}$ is an ${\rm{ORD}}$-length forcing iteration and  there is a proper class of $\mathbb{P}$-$\Sigma_{n+1}$-reflecting cardinals. Then,
\begin{enumerate}
\item $\mathcal{S}^{\Sigma_{n+1}}_{\mathbb{P}}\subseteq \mathcal{E}^{(n)}$  
\smallskip

%
%
\smallskip

\item If $\mathbb{P}$ is adequate and   $\Delta_{m+1}$-definable,  some $m\geq 1$,  then   
 $\mathcal{E}^{(m+n)}\subseteq\mathcal{S}^{\Sigma_{n+1}}_\mathbb{P}$.  
\end{enumerate}
In particular, if $\mathbb{P}$ is adequate   and $\Delta_{2}$-definable,      then  for every $n\geq 1$,   every $C^{(n+1)}$-extendible  is $\mathbb{P}$-$\Sigma_{n+1}$-supercompact.

\end{cor}

\section{Suitable iterations}
\label{suitit}

The following is a property enjoyed by many well-known ORD-length forcing iterations, such as Jensen's  canonical class forcing for obtaining the global GCH \cite{JensenForcing}, or the McAloon   class-forcing iteration for forcing $\mathrm{V=HOD}$ \cite{McAloon}.
The property will be needed to prove a general result (Theorem \ref{preservationCn}) about the preservation of $C^{(n)}$-extendibility.

\begin{defi}[Suitable iterations]\label{suitableiteration}
An ${\rm{ORD}}$-length forcing iteration $\mathbb{P}$   is \emph{suitable} if it is   the direct limit of an Easton support iteration\footnote{Recall that an Easton support iteration is a forcing iteration where direct limits are taken at inaccessible stages and inverse limits elsewhere.} $\langle \langle \mathbb{P}_\alpha:\,\alpha\in {\rm{ORD}}\rangle,\,\langle \dot{\mathbb{Q}}_\alpha:\,\alpha\in  {\rm{ORD}}\rangle\rangle$ with the property  that for each $\lambda\in {\rm{ORD}}$ 
there is some $\theta \in {\rm{ORD}}$ greater than $\lambda$ such that 
$$\Vdash_{\mathbb{P}_\nu}``\,\text{$\dot{\mathbb{Q}}_{\nu}$ is $\lambda$-directed closed}\,"$$
for all $\nu \geq \theta$.
\end{defi}


It is well-known that for suitable $\rm{ORD}$-length  forcing iterations $\mathbb{P}$, the forcing relation $\Vdash_{\mathbb{P}}$ is definable, the forcing theorem holds, and forcing with $\mathbb{P}$ preserves ZFC (see \cite{FRI} or \cite{Rie}).
The condition of eventual $\lambda$-directed closedness  in the definition above can be strengthened on a club proper class. Namely, 

\begin{prop}\label{clubclass}
Let  $\mathbb{P}$ be a  suitable iteration. The class $$C(\mathbb{P}):=\{\lambda: \forall\eta\geq\lambda, \Vdash_{\mathbb{P}_\eta}\text{$``\dot{\mathbb{Q}}_\eta$ is $\lambda$-directed closed''}\}$$
is a club class. 
\end{prop}
\begin{proof}
Closedness is obvious. As for unboundedness, fix any $\lambda$ and build inductively a sequence $\{\theta_n\}_{n\in\omega}$ of ordinals greater than $\lambda$ such that for all $\eta\geq\theta_{n+1}$, 
$$\Vdash_{\mathbb{P}_\eta}\text{`` $\dot{\mathbb{Q}}_\eta$ is $\theta_n$-directed closed''}.$$
Notice now that $\theta^*:=\sup_n\theta_n$ is an element of $C$.
\end{proof}

The next theorem establishes some sufficient conditions for the preservation of $C^{(n)}$-extendible cardinals under definable  iterations.
Recall that a partial ordering $\mathbb{P}$ is \emph{weakly homogeneous} if  for any $p,q\in \mathbb{P}$ there is an automorphism $\pi$ of $\mathbb{P}$ such that $\pi(p)$ and $q$ are compatible.  

In the case of definable ORD-length forcing iterations we define weak homogeneity as follows: 

\begin{defi}\label{WeakHomogeneityRevised}
A $\Gamma_m$-definable ORD-length forcing iteration $\mathbb{P}$ is  \emph{weakly homogeneous} if   there exists a  $\Gamma_m$ formula $\varphi(x,y, z_1, z_2, z_3)$  such that for every $\alpha$, $\mathbb{P}_\alpha$ forces that for every $\dot{p},\dot{q}\in \dot{\mathbb{Q}}_{\alpha}$, $\varphi(x,y, \dot{p},\dot{q},\alpha)$ defines an automorphism $\pi$ of $\dot{\mathbb{Q}}_{\alpha}$ such that $\pi(\dot{p})$ and $\dot{q}$ are compatible.  
\end{defi}

Let us note that if $\mathbb{P}$ is a weakly homogeneous $\Gamma_m$-definable iteration, and $\lambda$ is $\mathbb{P}$-$\Gamma_m$-reflecting, then $\mathbb{P}_\lambda$ is also a weakly homogeneous iteration in $V_\lambda$, with the same formula witnessing it. Also notice that, for every $\alpha$, $\mathbb{P}_\alpha$ forces that the remaining part of the iteration is weakly homogeneous. Indeed, in $V^{\mathbb{P}_\alpha}$, for every $p,q\in \mathbb{P}_{[\alpha, {\rm{ORD}})}$, the map $\pi$  given by: $\pi(x)=y$ iff for all $\beta >0$,
$$\Vdash_{\mathbb{P}_{ [ \alpha , \alpha +\beta)}}``\varphi (x(\beta), y(\beta), p(\beta), q(\beta), \beta)".$$
is a definable automorphism  of $\dot{\mathbb{P}}_{[\alpha, {\rm{ORD}})}$ (with parameters $\dot{p}$ and $\dot{q}$) such that $\mathds{1}_{\mathbb{P}_\alpha}$ forces $\pi(\dot{p})$ and $\dot{q}$ to be compatible.

\begin{theo}\label{preservationCn}
Let $m,n\geq 1$ and $m\leq n$. Let $\mathbb{P}$ be an adequate  $\Delta_{m+1}$-definable  and weakly homogeneous suitable  iteration. Suppose   there is a proper class of   $\mathbb{P}$-$\Sigma_{n+1}$-reflecting cardinals.  
If $\delta$ is a $\mathbb{P}$-$\Sigma_{n+1}$-supercompact cardinal, then 
$$\Vdash_{\mathbb{P}} \text{$``\delta$ is \Cn-extendible''}.$$ 
\end{theo}

\begin{proof}
 Suppose $G$ is $\mathbb{P}$-generic over $\rm{V}$. By Corollary \ref{equivalence} and Lemma \ref{PsigmareflectingCncardinal}, it is sufficient to take an arbitrary   $\mathbb{P}$-$\Sigma_{n+1}$-reflecting cardinal $\lambda>\delta$  
 and any $\alpha <\lambda$,  and find a $\mathbb{P}$-$\Sigma_{n+1}$-reflecting cardinal $\bar{\lambda}$,  
 ordinals $\bar{\delta}, \bar{\alpha}<\bar{\lambda}$, and an elementary embedding $j: V[G]_{\bar{\lambda}}\longrightarrow V[G]_\lambda$ such that $\crit(j)=\bar{\delta}$, $j(\bar{\delta})=\delta$, and $j(\bar{\alpha})=\alpha$. 

So pick a  $\mathbb{P}$-$\Sigma_{n+1}$-reflecting cardinal $\lambda>\delta$  
and   any $\alpha <\lambda$.  
Since $\delta$ is  $\mathbb{P}$-$\Sigma_{n+1}$-supercompact there exist $\bar{\delta}<\bar{\lambda}< \delta$ and $\bar{\alpha}<\bar{\lambda}$, together with an elementary embedding $j: V_{\bar{\lambda}}\longrightarrow V_{\lambda}$ such that 
\begin{itemize}
\item[(i)] $\crit(j)=\bar{\delta}$ and  $j(\bar{\delta})=\delta$ 
\item[(ii)] $j(\bar{\alpha})=\alpha$ 
\item[(iii)] $\bar{\lambda}$ is $\mathbb{P}$-$\Sigma_{n+1}$-reflecting.
\end{itemize}
It will then suffice to show that $j$ can be lifted to  an elementary embedding $j:V_{\bar{\lambda}}[G_{\bar{\lambda}}]\longrightarrow V_\lambda[G_\lambda]$, for then, since both $\lambda$ and $\bar{\lambda}$ are $\mathbb{P}$-reflecting, we have that $V_\lambda[G_\lambda]=V[G]_\lambda$ and  $V_{\bar{\lambda}}[G_{\bar{\lambda}}]=V[G]_{\bar{\lambda}}$. 

\smallskip

The iterations $\mathbb{P}_{\bar{\lambda}}$ and $\mathbb{P}_{{\lambda}}$ factorize as follows:
\begin{itemize}
\item[(i)] $\mathbb{P}_{\bar{\lambda}}\cong \mathbb{P}_{\bar{\delta}}\ast \mathbb{Q}$ with $|\mathbb{Q}|\leq \bar{\lambda}$.
\item[(ii)]  $\mathbb{P}_{{\lambda}}\cong \mathbb{P}_{{\delta}}\ast \mathbb{Q}^*$ with 
$$\Vdash_{\mathbb{P}_\delta} \text{``$\mathbb{Q}^*$ is weakly homogeneous and $\delta$-directed closed''}.$$
\end{itemize}
Indeed, (i) is clear since $\bar{\lambda}$ is a $\beth$-fixed point. For (ii), since  $\mathbb{P}$ is weakly homogeneous and $\lambda$ is $\mathbb{P}$-$\Sigma_{n+1}$-reflecting, $\mathbb{P}_\lambda$ is weakly homogeneous in $V_\lambda$, and therefore $\mathbb{P}_\delta$ forces that $\mathbb{Q}^\ast$ is weakly homogeneous. Thus, we only need to see that 
 $\Vdash_{\mathbb{P}_\delta} ``\text{$\mathbb{Q}^*$ is $\delta$-directed closed''}$. 
 
 Recall from Proposition \ref{clubclass} that the class $$C(\mathbb{P}):=\{\mu : \forall\eta\geq\mu,\, \Vdash_{\mathbb{P}_\eta}\text{`` $\dot{\mathbb{Q}}_\eta$ is $\mu$-directed closed''}\}$$
is a club class. Thus, it will be sufficient to show that $\delta$ is a limit point of $C(\mathbb{P})$, and therefore it belongs to $C(\mathbb{P})$.  So, let $\mu<\delta$ and notice that since $\mathbb{P}$ is a suitable iteration, the sentence $\varphi(\mu)$ asserting: 
$$\exists \theta >\mu\, \forall\eta\geq\theta \; ( \; \Vdash_{\mathbb{P}_\eta}\text{``$\dot{\mathbb{Q}}_\eta$ is $\mu$-directed closed''})$$
holds in $\rm{V}$. Since $\mathbb{P}$ is $\Delta_{m+1}$-definable,  
$\varphi(\mu)$ is easily seen to be equivalent to the $\Sigma_{m+2}$ sentence
$$\exists \theta>\mu \, \forall\eta\geq \theta\, \forall \alpha >\eta \, (\alpha \in C^{(m)}\to$$ $$   V_\alpha\vDash ``\Vdash_{\mathbb{P}_\eta}``\,\text{$\dot{\mathbb{Q}}_\eta$ is $\mu$-directed closed}\, " ").$$
Since $\delta$ is  a $\Sigma_{n+2}$-correct cardinal (by Lemma \ref{SigmanSupercompIsCorrect} and Proposition \ref{PropertiesPsigmasupercom}), and $m\leq n$, there must be a witness for $\varphi(\mu)$ below $\delta$. Arguing inductively, we define an $\omega$-sequence of ordinals above $\mu$ with limit in $C(\mathbb{P})$.  This shows that $C(\mathbb{P})$ is unbounded in $\delta$, as wanted.

\medskip

Since $\delta\in C^{(n+2)}$, and  $j$ is elementary with $j(\bar{\delta})=\delta$, we have that $j(\mathbb{P}_{\bar{\delta}})=\mathbb{P}_\delta$. Also, since $\bar{\delta}$ is  the critical point of $j$, we have that $j``G_{\bar{\delta}}=G_{\bar{\delta}} \subseteq G_\delta$, and so  $j\restriction V_{\bar{\lambda}}$ can be lifted to an elementary embedding
$$j:V_{\bar{\lambda}}[G_{\bar{\delta}}]\longrightarrow V_{\lambda}[G_\delta].$$

Let us  denote by $G_{[\bar{\delta},\bar{\lambda})}$ and $G_{[{\delta},{\lambda})}$ the filters $G\cap \mathbb{Q}$ and $G\cap \mathbb{Q}^*$, respectively. Notice that these filters are generic  for $\mathbb{Q}$ and $\mathbb{Q}^*$ over $V_{\bar{\lambda}}[G_{\bar{\delta}}]$ and $V_{\lambda}[G_\delta]$, respectively. In order to lift the embedding $j$ to  the further generic extension $V_{\bar{\lambda}}[G_{\bar{\lambda}}]=V_{\bar{\lambda}}[G_{\bar{\delta}}][G_{[\bar{\delta},\bar{\lambda})}]$, notice first that $j``G_{[\bar{\delta},\bar{\lambda})}$ is a directed subset of $\mathbb{Q}^\ast$  of cardinality  $\leq \bar{\lambda}$. Also, since $j\upharpoonright V_{\bar{\lambda}}\in V_\lambda[G_\delta]$, and $j``G_{[\bar{\delta},\bar{\lambda})}$ can be computed from $j\upharpoonright V_{\bar{\lambda}}$ and $G_{[\bar{\delta},\bar{\lambda})}$, we have that $j``G_{[\bar{\delta},\bar{\lambda})}\in V_\lambda[G_\delta]$.  
Therefore, since $\mathbb{Q}^*$ is a $\delta$-directed closed forcing notion in $V_{\lambda}[G_\delta]$, there is some condition $p\in\mathbb{Q}^*$ such that $p\leq q$, for every $q\in j``G_{[\bar{\delta},\bar{\lambda})}$. Thus, $p$ is a master condition in $\mathbb{Q}^*$ for the embedding $j$ and the generic filter $G_{[\bar{\delta},\bar{\lambda})}$. So, if $H\subseteq \mathbb{Q}^*$ is a generic filter over $V_{\lambda}[G_\delta]$ containing $p$, then  $j$ can be lifted to an elementary embedding
$$j: V_{\bar{\lambda}}[G_{\bar{\lambda}}]\longrightarrow V_{\lambda}[G_\delta\ast H].$$
\begin{claim}
In $V[G]$ there exists some generic filter $H\subseteq \mathbb{Q}^*$ over $V_{\lambda}[G_\delta]$ containing $p$ such that $V_{\lambda}[G_\delta\ast H]=V_{\lambda}[G_\lambda]$.
\end{claim}
\begin{proof}[Proof of claim]
By (ii) above,  $\mathbb{Q}^*$ is a weakly homogeneous class forcing in $V_{\lambda}[G_\delta]$. Thus, the set of conditions $r\in \mathbb{Q}^*$ for which there is an automorphism $\pi$ of $\mathbb{Q}^*$ that is definable in $V_\lambda [G_\delta]$ and such that $\pi(r)\leq p$ is dense.  Pick such an $r$ in $G_{[\delta,\lambda)}$ and such an automorphism $\pi$. Now, notice that the filter $H$ generated by the set $\pi``G_{[\delta,\lambda)}$ contains $\pi(r)$ and therefore it contains $p$. Since $H$ is definable by means of $\pi$ and $G_{[\delta,\lambda)}$, and also $\pi$ is definable in $V_\lambda[G_\delta]$,  we conclude that $V_{\lambda}[G_\delta\ast H]=V_{\lambda}[G_\lambda]$.
\end{proof}

By taking $H\subseteq \mathbb{Q}^*$ as in the claim above, we thus obtain a lifting 
$$j:V_{\bar{\lambda}}[G_{\bar{\lambda}}]\longrightarrow V_{\lambda}[G_\lambda]$$
as wanted.
\end{proof}

\begin{remark}\label{RemarkOnLocality}
Note that the above proof shows more than what  is stated in Theorem~\ref{preservationCn}. Specifically,
what is proved is the following local result:  every  elementary embedding (in $V$) witnessing some partial degree of $\mathbb{P}$-$\Sigma_{n+1}$-supercompactness lifts to an elementary embedding (in $V^{\mathbb{P}}$) witnessing the same partial degree of $C^{(n)}$-extendibility. 
\end{remark}

\begin{cor} 
\label{maincor}
Suppose $n\geq 1$. Let $\mathbb{P}$ be an adequate  $\Delta_2$-definable and  weakly homogeneous suitable iteration. 
Suppose 
there is  a proper class of   $\mathbb{P}$-$\Sigma_{n+1}$-reflecting cardinals. 
If $\delta$ is a $C^{(n+1)}$-extendible cardinal, then $$\Vdash_{\mathbb{P}}``\, \text{$\delta$ is \Cn-extendible}\,".$$
\end{cor}
\begin{proof}
Since $\mathbb{P}$ is adequate,  $\Delta_2$-definable, and there exists a proper class of $\mathbb{P}$-$\Sigma_{n+1}$-reflecting cardinals, Corollary \ref{coroscext} implies that $\delta$ is $\mathbb{P}$-$\Sigma_{n+1}$-supercompact. Now, Theorem \ref{preservationCn} applies to get the desired result. 
\end{proof}

\bigskip

Let $\mathbb{P}$ be an adequate    $\Delta_{m+1}$-definable weakly homogeneous 
suitable iteration.
Let us  look next into the conditions under which,  
for some $n\geq 1$,  there exists  a  proper class of  $\mathbb{P}$-$\Sigma_{n+1}$-reflecting cardinals (this was one of the assumptions  of Theorem \ref{preservationCn}). 

First, the class $C^{(n)}_{\mathbb{P}}$ is closed, unbounded, and  $\Pi_{m+n}$-definable (Proposition \ref{complexityclubclass}). Also, the class $D$ of ordinals $\alpha$ such that $\mathbb{P}_\alpha=\mathbb{P}\cap V_\alpha$ is   
 $\Pi_{m+1}$-definable  (cf. the proof of Proposition \ref{ComplexityPsigmareflecting}). Observe that $D$ is a proper class by assumption, hence 
$\overline{D}$ (the closure of 
$D$)  is a club proper class.  Further,   the class $K$ of $\mathbb{P}$-reflecting cardinals 
is unbounded.

Indeed, fix any ordinal $\alpha$ and define inductively a sequence  of ordinals $\langle \alpha_n :n<\omega\rangle$  in $D$, with $\alpha_0=\alpha$ and $\forces_\mathbb{P}``V[\dot{G}]_{\alpha_n}\s V_{\alpha_{n+1}}[\dot{G}_{\alpha_{n+1}}]$'',  as follows: Suppose that $\alpha_n$ has already been  defined.  
Since $\mathbb{P}$ is suitable,  it preserves $\mathrm{ZFC}$, and therefore $\mathbb{P}$ forces that $V[\dot{G}]_{\alpha_n}$ is a set. Suppose for a moment that there is $\mathbb{P}$-name $\tau$ such that $\mathbb{P}$ forces $``\tau=V[\dot{G}]_{\alpha_n}$''. Let $\beta$ be an ordinal greater than the rank of $\tau$ and such that $p\in\mathbb{P}_\beta$ and $\mathbb{P}_\beta=\mathbb{P}\cap V_\beta$.  Then,  $\forces_{\mathbb{P}}\tau\in V_\beta[\dot{G}_\beta]$ and thus $\forces_{\mathbb{P}} ``V[\dot{G}]_{\alpha_n}\s V_{\alpha_{n+1}}[\dot{G}_{\alpha_{n+1}}]$'', where  $\beta:=\alpha_{n+1}$. So let us show that such a  $\mathbb{P}$-name $\tau$ exists.  
 Let $\beta$ be the first ordinal such that $\mathbb{P}_\beta=\mathbb{P}\cap V_\beta$ and  there exist  $p\in\mathbb{P}_\beta$ and  $\tau_p\in V_\beta^{\mathbb{P}_\beta}$ such that $p\forces_{\mathbb{P}}``\tau_p=V[\dot{G}]_{\alpha_n}$''. Let now $q\in \mathbb{P}_\beta$ be arbitrary. By weak homogeneity of $\mathbb{P}$  let $\pi$ be an automorphism  of $\mathbb{P}$ such that $\pi(p)$ and $q$ are compatible. Note that $\pi(\tau_p)$ is a $\mathbb{P}_\beta$-name (cf. Definition~\ref{WeakHomogeneityRevised}) and that $\pi(p)\forces_{\mathbb{P}}``\pi(\tau_p)=V[\dot{G}]_{\alpha_n}"$.  Now, letting $r\leq_{\mathbb{P}_\beta} \pi(p), q$  we have that $r\forces_{\mathbb{P}_\beta}``\pi(\tau_p)=V[\dot{G}]_{\alpha_n}$''. This shows that the set of conditions $q\in \mathbb{P}_\beta$ for which there is a $\mathbb{P}_\beta$-name $\tau_q$  such that $q\forces_{\mathbb{P}_\beta}``\tau_q=V[\dot{G}]_{\alpha_n}$'' is dense in $\mathbb{P}_\beta$. So, let $A$ be a maximal antichain of $\mathbb{P}_\beta$ contained in this dense set.  Now  we may combine all the names $\tau_q$, for $q\in A$, into a single $\mathbb{P}_\beta$-name  $\tau$ so that $q\forces_{\mathbb{P}_\beta}``\tau=V[\dot{G}]_{\alpha_n}"$, for all $q\in A$. Hence, $\forces_{\mathbb{P}_\beta}``\tau=V[\dot{G}]_{\alpha_n}"$, and therefore $\forces_{\mathbb{P}}``\tau=V[\dot{G}]_{\alpha_n}"$.


Finally, let $\lambda$ be the supremum of the $\alpha_n$'s. 
In this case we have that $\Vdash_{\mathbb{P}}``V[\dot{G}]_\lambda \subseteq V_\lambda [\dot{G}_\lambda]"$, and so $\lambda$ is $\mathbb{P}$-reflecting (see Remark \ref{remarkPReflecting}). Hence, $K$ is a club proper class (Proposition \ref{Preflectingclosed}). Also, by virtue of Proposition \ref{ComplexityPreflecting}, $K$ is   $\Pi_{m+2}$-definable.

 For each $n\geq 1$, set $\mathcal{K}_{n}:=C^{(n)}_{\mathbb{P}} \cap {K}\cap \overline{D}$.

\begin{prop}\label{complexityofKn} 
Let $\mathbb{P}$ be an adequate $\Delta_{m+1}$-definable weakly homogeneous suitable iteration.
Then $\mathcal{K}_{n+1}$ (resp. $\mathcal{K}_1$) is a $\Pi_{m+n+1}$-definable (resp. $\Pi_{m+2}$) club class such that $\mathcal{K}_{n+1}\cap \mathrm{Reg}$ (resp. $\mathcal{K}_1\cap \rm{Reg}$) is contained in the class of  $\mathbb{P}$-$\Sigma_{n+1}$-reflecting cardinals (resp. $\mathbb{P}$-$\Sigma_{1}$-reflecting cardinals). 
\end{prop}
\begin{proof}
Clearly, $\mathcal{K}_{n+1}$ is a club class. Also, by the above discussion, $\mathcal{K}_{n+1}$ is $\Pi_{m+n+1}$, and $\mathcal{K}_1$ is $\Pi_{m+2}$.  Let  $\kappa\in \mathcal{K}_{n+1}\cap \mathrm{Reg}$. To show that $\kappa$ is $\mathbb{P}$-$\Sigma_{n+1}$-reflecting it is enough to verify that $\kappa\in D$.

Since $\kappa$ is an accumulation point of $D$, $\mathbb{P}\cap V_\kappa\subseteq \mathbb{P}_\kappa$. Conversely, 
 $\kappa$ is inaccessible, hence  $\mathbb{P}_\kappa$ is the direct limit of the previous stages, and so $\mathbb{P}_\kappa\subseteq \mathbb{P}\cap V_\kappa$. Hence $\kappa\in D$. The proof for $\mathcal{K}_1$ is the same. 
\end{proof}

Recall that if $\kappa$ is a cardinal and $\alpha$   an ordinal, 
\begin{itemize}
\item $\kappa$ is $0$-Mahlo if $\kappa$ is inaccessible;
\item  $\kappa$ is $(\alpha+1)$-Mahlo iff $\{\lambda<\kappa\mid \text{$\lambda$ is $\alpha$-Mahlo}\}$ is stationary;
\item In case $\alpha>0$ is a limit ordinal, $\kappa$ is $\alpha$-Mahlo iff $\kappa$ is $\beta$-Mahlo for every $\beta<\alpha$.
\end{itemize}

\begin{defi}
For $n<\omega$ and $\alpha\in\rm{ORD}$, the class   \emph{${\rm{ORD}}$ is $\mathbf{\Gamma_n}$-$(\alpha +1)$-Mahlo} if every $\mathbf{\Gamma_n}$-definable\footnote{i.e., $\Gamma_n$-definable with parameters (see Notation~\ref{NotationGamma}).} club proper class of ordinals contains an $\alpha$-Mahlo cardinal.
\end{defi}

Note that a cardinal is  Mahlo if and only if it is $1$-Mahlo, hence $\mathrm{ORD}$ is $\mathbf{\Gamma_n}$-Mahlo if and only if $\mathrm{ORD}$ is $\mathbf{\Gamma_n}$-1-Mahlo. 

\begin{defi}
A proper class $S$ of ordinals is \emph{$\mathbf{\Gamma_n}$-stationary} if $S$ intersects every $\mathbf{\Gamma_n}$-definable club proper class of ordinals.
\end{defi}

%

\begin{prop}\label{propMahlo}
Let $\mathbb{P}$ be an adequate $\Delta_{m+1}$-definable weakly homogeneous suitable iteration. If $\rm{ORD}$ is   $\mathbf{\Pi_{m+n+1}}$-Mahlo, then  the class of $\mathbb{P}$-$\Sigma_{n+1}$-reflecting cardinals is  $\mathbf{\Pi_{m+n+1}}$-stationary, so it is a proper class.
\end{prop}
\begin{proof}
Let $\mathcal{C}$ be a  $\mathbf{\Pi_{m+n+1}}$ club proper class of ordinals. By Proposition \ref{complexityofKn}, $\mathcal{C}\cap \mathcal{K}_{n+1}$ is a  $\mathbf{\Pi_{m+n+1}}$ club proper class such that $\mathcal{C}\cap \mathcal{K}_{n+1}\cap \rm{Reg}$ is contained in the class of $\mathbb{P}$-$\Sigma_{n+1}$-reflecting cardinals. Thus, the class of $\mathbb{P}$-$\Sigma_{n+1}$-reflecting cardinals is  $\mathbf{\Pi_{m+n+1}}$-stationary.
To show the class is unbounded, given any cardinal $\kappa$, let $\mathcal{C}:=C^{(n+m)}\setminus \kappa^+$. Then $\mathcal{C}$ is a $\mathbf{\Pi_{m+n+1}}$-definable club proper class  which, as before, contains a $\mathbb{P}$-$\Sigma_{n+1}$-reflecting cardinal.  \end{proof}

\begin{lemma}\label{LemmaAdequate}
Let $\mathbb{P}$ be a $\Delta_{m+1}$-definable suitable iteration. If $\rm{ORD}$ is $\mathbf{\Pi_{m+1}}$-Mahlo, then $\mathbb{P}$ is adequate.
\end{lemma}
\begin{proof}
Since $\mathbb{P}$ is $\Delta_{m+1}$-definable, $\mathbb{P}^{V_\kappa}=\mathbb{P}\cap V_\kappa$ for every $\kappa \in C^{(m+1)}$. Since $\rm{ORD}$ is $\mathbf{\Pi_{m+1}}$-Mahlo, and direct limits are taken at inaccessible stages of the iteration, the class  of inaccessible cardinals $\kappa$ such that $\mathbb{P}^{V_\kappa}=\mathbb{P}_\kappa$ is a proper class. Thus,  $\mathbb{P}$ is adequate.
\end{proof}

The following corollaries now follow immediately from Theorem \ref{preservationCn}, Corollary \ref{maincor},  Proposition \ref{propMahlo} and Lemma \ref{LemmaAdequate}:

\begin{cor}
\label{coroMahlo}
Let $m,n\geq 1$ and $m\leq n$. Let $\mathbb{P}$ be 
$\Delta_{m+1}$-definable  weakly homogeneous suitable  iteration. Suppose   ${\rm{ORD}}$ is $\mathbf{\Pi_{m+n+1}}$-Mahlo.
If $\delta$ is a $\mathbb{P}$-$\Sigma_{n+1}$-supercompact cardinal, then 
$$\Vdash_{\mathbb{P}} \text{$``\delta$ is \Cn-extendible''}.$$ 
\end{cor}
 
\begin{cor} 
\label{maincor2}
Let $n\geq 1$ and let $\mathbb{P}$ be 
$\Delta_2$-definable weakly homogeneous suitable iteration. 
Suppose   ${\rm{ORD}}$ is $\mathbf{\Pi_{n+2}}$-Mahlo.
If $\delta$ is a $C^{(n+1)}$-extendible cardinal, then $$\Vdash_{\mathbb{P}}``\, \text{$\delta$ is \Cn-extendible}\,".$$
\end{cor}

\section{Vop\v{e}nka's  principle and suitable iterations}
Let us recall the following characterizations of Vop\v{e}nka's Principle, as well as of its restriction to definable classes of structures of a given complexity, in terms of   $C^{(n)}$-extendible cardinals. 

\begin{theo}[\cite{BagEtAl}]\label{VopenkaBagaria}
The following are equivalent:
\begin{enumerate}
\item $\rm{VP}$
\item For every  $n\geq 1$ there exists a \Cn-extendible cardinal.
\end{enumerate}
\end{theo}

\begin{theo}[\cite{Bag}]\label{VopenkaBagaria}
Let $n\geq 1$. The following are equivalent:
\begin{enumerate}
\item $\VP(\mathbf{\Pi_{n+1}})$, i.e., $\VP$ restricted to proper classes of structures that are ${\Pi_{n+1}}$-definable with parameters.
\item There exists a proper class of \Cn-extendible cardinals.
\end{enumerate}
\end{theo}

For $n=0$, the  equivalence  is between $\VP (\mathbf{\Pi_1})$ and the existence of a proper class of supercompact cardinal. The \emph{lightface} version also holds. Namely, $\VP(\Pi_{n+1})$ is equivalent to the existence of a  \Cn-extendible cardinal (see \cite{Bag}).

\medskip

Vop\v{e}nka's Principle can be also characterized in terms of the existence of  $\mathbb{P}$-$\Sigma_n$-supercompact cardinals (cf. Theorem \ref{equivalencesvopenka} and Theorem \ref{VopParam}). The following lemma will be useful for this purpose.

\begin{lemma}\label{OrdMahlo} 
Let $n\geq 1$. Then 
$\VP(\mathbf{\Pi_n})$ implies   
that  $\mathrm{ORD}$ is $\mathbf{\Sigma_{n+1}}$-$(\kappa+1)$-Mahlo, for every cardinal $\kappa$.\footnote{The lightface version of the statement also holds.}
\end{lemma}
\begin{proof}
Let us  prove the lemma for $n>1$. The case $n=1$  is similar, using the fact that $\VP(\mathbf{\Pi_1})$ is equivalent to the existence of a proper class of supercompact cardinal, and that every supercompact cardinal belogs to $C^{(2)}$. So, let $n>1$ and assume that $\VP(\mathbf{\Pi_{n}})$ holds.  

Let $\mathcal{C}$ be a $\mathbf{\Sigma_{n+1}}$-definable club proper class of ordinals and let $\varphi(x,\vec{a})$ be some $\Sigma_{n+1}$-formula defining it. Let $\kappa$ be a cardinal and  $\lambda$ be a $C^{(n-1)}$-extendible cardinal with $\vec{a}\in V_\lambda$ and $\kappa
<\lambda$. Note that  this $\lambda$ exists by virtue of Theorem \ref{VopenkaBagaria}. We claim that $\mathcal{C}\cap \lambda$ is unbounded in $\lambda$. For if $\alpha <\lambda$, then the sentence ``$\exists \beta >\alpha (\beta \in \mathcal{C})$" is $\Sigma_{n+1}$ (with $\vec{a}$ as parameters), hence it is true in $V_\lambda$ because $\lambda$ is $C^{(n-1)}$-extendible and so it belongs to $C^{(n+1)}$ (\cite[Proposition 3.4]{Bag}). Since $\mathcal{C}$ is closed,  $\lambda\in \mathcal{C}$. Since $\lambda$ is $\lambda$-Mahlo, hence also $\kappa$-Mahlo, the result follows.
\end{proof}

\begin{theo}\label{equivalencesvopenka}
The following are equivalent:
\begin{enumerate}
\item $\VP$ holds.
\item For every $n\geq 1$ and every  
definable weakly homogeneous suitable iteration $\mathbb{P}$, there is a proper class of $\mathbb{P}$-$\Sigma_{n}$-reflecting cardinals   and there exists  a  proper class of $\mathbb{P}$-$\Sigma_n$-supercompact cardinal.
\end{enumerate}
\end{theo}
\begin{proof}
$(1) \Rightarrow (2)$: Let $n\geq 1$ and let $\mathbb{P}$ be  $\Delta_{m+1}$-definable weakly homogeneous suitable iteration, some $m\geq 1$.   By Lemma \ref{LemmaAdequate}, $\mathbb{P}$ is adequate.  
By  Lemma \ref{OrdMahlo},   
 ${\rm{ORD}}$ is $\mathbf{\Pi_{m+n+1}}$-Mahlo and so  Proposition \ref{propMahlo} yields a proper class of $\mathbb{P}$-$\Sigma_{n+1}$-reflecting cardinals. Also, by  Corollary \ref{coroscext} and Theorem \ref{VopenkaBagaria}   we infer that there is a proper class of $\mathbb{P}$-$\Sigma_n$-supercompact cardinal.

$(2)\Rightarrow (1)$: Since (2) holds, by  Corollary \ref{coroscext} there exists a proper class of  $C^{(n)}$-extendible cardinals, for every $n\geq 1$, hence by Theorem \ref{VopenkaBagaria}, $\rm{VP}$ holds.
\end{proof}

The following gives more precise information about the relationship between $\mathbb{P}$-$\Sigma_n$-supercompact cardinals and fragments of $\VP$.

\begin{theo}\label{VopParam}
Let  $m,n\geq 1$ and let $\mathbb{P}$ be an ${\rm{ORD}}$-length forcing  iteration. Then, 
\begin{enumerate}
\item Assume there is a proper class of $\mathbb{P}$-$\Sigma_{n+1}$-reflecting cardinals. If there is a $\mathbb{P}$-$\Sigma_{n+1}$-supercompact cardinal, then $\VP(\Pi_{n+1})$ holds; and if there is a proper class of $\mathbb{P}$-$\Sigma_{n+1}$-supercompact cardinals, then $\VP(\mathbf{\Pi_{n+1}})$ holds.
\item  If $\mathbb{P}$ is  
$\Delta_{m+1}$-definable,  some $m\geq 1$,  suitable and weakly homogeneous, then 
$\VP(\mathbf{\Pi_{m+n+1}})$ implies  the existence of  a proper class of  $\mathbb{P}$-$\Sigma_{n+1}$-reflecting cardinals and a proper class of $\mathbb{P}$-$\Sigma_{n+1}$-supercompact cardinals.
\end{enumerate}
\end{theo}
\begin{proof}
Item (1) is a direct consequence of Corollary \ref{coroscext} and our remarks following Theorem \ref{VopenkaBagaria}. As for (2),  
Theorem \ref{VopenkaBagaria} show that $\VP(\mathbf{\Pi_{m+n+1}})$ 
is equivalent to the existence of a proper class of $C^{(m+n)}$-extendible cardinals.   
 Also, Lemma \ref{OrdMahlo}   shows that $\VP(\mathbf{\Pi_{m+n+1}})$ 
 implies that ORD is $\mathbf{\Sigma_{m+n+2}}$-Mahlo, hence Lemma \ref{LemmaAdequate} implies that $\mathbb{P}$ is adequate. 
  Thus, by Proposition \ref{propMahlo},
 there exists a proper class of $\mathbb{P}$-$\Sigma_{n+1}$-reflecting 
 cardinals. By Corollary \ref{coroscext} 
 there exists also a proper class of $\mathbb{P}$-$\Sigma_{n+1}$-supercompact cardinals. 
 \end{proof}

We end this section by proving a level-by-level version of Brooke-Taylor's result on the preservation of Vop\v{e}nka's Principle under definable suitable iterations.

\begin{theo}\label{AlmostVopenka}
Let $m,n\geq 1$ and $m\leq n$. Let $\mathbb{P}$ be a 
$\Delta_{m+1}$-definable weakly homogeneous suitable iteration.   If   $\VP(\mathbf{\Pi_{m+n+1}})$ holds, then $$\Vdash_\mathbb{P}\text{$``\mathrm{VP}(\mathbf{\Pi_{n+1}})$ holds''}.$$
\end{theo}

\begin{proof}
By Theorem \ref{VopParam}(2), there is a proper class of $\mathbb{P}$-$\Sigma_{n+1}$-reflecting cardinals, as well as a proper class of $\mathbb{P}$-$\Sigma_{n+1}$-supercompact cardinals. Also, by Lemma \ref{LemmaAdequate}, $\mathbb{P}$ is adequate. Now the result follows combining Theorem \ref{preservationCn} and the remarks just after Theorem \ref{VopenkaBagaria}. 
\end{proof}

\begin{cor}[\cite{Broo}]\label{VopenkaIndes}
Let   $\mathbb{P}$ be a   definable weakly homogeneous suitable iteration. If $\rm{VP}$ holds in $\rm{V}$, then $\rm{VP}$ holds in $V^{\mathbb{P}}$.
\end{cor}

Our version of Brooke-Taylor's result differs from the original one in that we require the weak homogeneity of $\mathbb{P}$. 
However,  our proof shows more than Brooke-Taylor's, for it shows that \emph{every} relevant elementary embedding from the ground model lifts to an elementary embedding in the forcing extension (recall Remark~\ref{RemarkOnLocality}).
Even though weak homogeneity  holds for a wide family of forcing notions, it puts some restrictions on the sort of statements that can be forced. One example is   $``\rm{V}=\rm{HOD}"$. 
 In Section 12 we will address this problem and will 
 prove Theorem \ref{AlmostVopenka} without the weak homogeneity assumption (Theorem \ref{VopenkaNonHom}).
We are very grateful to Andrew Brooke-Taylor for his valuable  comments on this matter. 

\section{\Cn-extendible cardinals and fitting iterations}\label{SectionFitting}
While our main corollary from Section \ref{suitit} (cf. Corollary \ref{maincor2}) applies to a wide class of forcing iterations, it requires that the cardinal we begin with is $C^{(n+1)}$-extendible. In this section we argue that this extra assumption can be avoided by restricting  to a certain subclass of forcing iterations, which  
nevertheless are still general enough. 

\begin{defi}\label{NewFitting}
 An $\rm{ORD}$-length 
 forcing iteration $\mathbb{P}$ is \emph{fitting} if it is $\Delta_2$-definable, suitable and weakly homogeneous, and there is a $\Delta_2$-definable  $\mathbf{\Sigma_{2}}$-stationary  class $K_\mathbb{P}$   of regular  $\mathbb{P}$-reflecting cardinals $\kappa$ such that $\mathbb{P}\cap V_\kappa=\mathbb{P}_\kappa$.
\end{defi}

\begin{remark}\label{RemarkonKP}
Note that every fitting iteration is adequate.
Also notice that 
$K_\mathbb{P}\cap C^{(1)}$ is a $\Delta_2$-definable proper class of inaccessible cardinals.
\end{remark}

\begin{theo}\label{2fittingsupercompact}
Suppose that $\mathbb{P}$ is a fitting iteration.  If $\delta$ is a supercompact cardinal  in $C(\mathbb{P})$, then $\forces_\mathbb{P}\text{$``\delta$ is supercompact''}.$ \footnote{For the definition of $C(\mathbb{P})$, see Proposition \ref{clubclass}.}
\end{theo}
\begin{proof}
Let $\delta$ be a supercompact  cardinal in $C(\mathbb{P})$ and suppose $G$ is a $\mathbb{P}$-generic filter over $\rm{V}$. 
It will be sufficient to  verify that for each $\lambda\in K_\mathbb{P}\cap C^{(1)}$   there are $\bar{\lambda}\in K_\mathbb{P}\cap C^{(1)}$  and $\bar{\delta}$ with $\bar{\delta}< \bar{\lambda}<\delta,$ and an elementary embedding
$$j\colon V_{\bar{\lambda}}[G_{\bar{\lambda}}]\rightarrow V_\lambda[G_\lambda]$$
with $\crit(j)=\bar{\delta}$   and  $j(\bar{\delta})=\delta$. For since $\lambda $ and $\bar{\lambda}$ are $\mathbb{P}$-reflecting  we then have $V_{\bar{\lambda}}[G_{\bar{\lambda}}]=V[G]_{\bar{\lambda}}$ and $V_{\lambda}[G_{\lambda}]=V[G]_{\lambda}$. 

So let $\lambda$ be as above. Since $\delta$ is $\Sigma_1$-supercompact, given  $\mu\in C^{(1)}$ with $\mu>\lambda$, there is $\bar{\mu}\in C^{(1)}$ and $\bar{\delta}<\bar{\lambda}<\bar{\mu}$, and an elementary embedding $j:V_{\bar{\mu}}\rightarrow V_\mu$ with $\crit(j)=\bar{\delta}$, $j(\bar{\delta})=\delta$ and $j(\bar{\lambda})=\lambda$.

Since $K_\mathbb{P}$ is $\Delta_2$-definable, $V_\mu \models \lambda\in K_\mathbb{P}\cap C^{(1)}$. Hence, by elementarity and $\Sigma_1$-correctness of $\bar{\mu}$, $\bar{\lambda}\in K_\mathbb{P}\cap C^{(1)}$.

We only need to show how to lift   $j\restriction V_{\bar{\lambda}}:V_{\bar{\lambda}}\rightarrow V_\lambda$ to an elementary embedding $V_{\bar{\lambda}}[G_{\bar{\lambda}}]\rightarrow V_\lambda[G_\lambda].$ Since $\mathbb{P}$ is $\Delta_2$-definable,  $\mathbb{P}^{V_\lambda}=\mathbb{P}\cap V_\lambda=\mathbb{P}_\lambda$, where the rightmost equality follows from $\lambda\in K_\mathbb{P}$. Analogously, the same  holds for $\mathbb{P}^{V_{\bar{\lambda}}}$. Thus, we have:
\begin{itemize}
\item[(i)] $\mathbb{P}_{\bar{\lambda}}\cong \mathbb{P}_{\bar{\delta}}\ast \dot{\mathbb{Q}}$, with $|\dot{\mathbb{Q}}|=\bar{\lambda}$.
\item[(ii)] $\mathbb{P}_\lambda\cong \mathbb{P}_{\delta}\ast \dot{\mathbb{Q}}^*$ and $\forces_{\mathbb{P}_{\delta}}\text{$``\dot{\mathbb{Q}}^*$ is weakly homogeneous''}$.
\end{itemize}
Since  $\delta\in C(\mathbb{P})$,  $\forces_{\mathbb{P}_\delta}\text{$``\dot{\mathbb{Q}}^*$ is $\delta$-directed closed''}.$
We may now proceed as  in  the proof of Theorem  \ref{preservationCn} to get the desired elementary embedding
$j:V_{\bar{\lambda}}[G_{\bar{\lambda}}]\rightarrow V_\lambda[G_\lambda].$
\end{proof}

\begin{lemma}\label{2fittingextendibles}
Suppose that $\mathbb{P}$ is a fitting iteration and $\delta$ is an extendible cardinal. Then, $\forces_\mathbb{P}\text{$``\delta$ is extendible''}.$
\end{lemma}

\begin{proof}
Suppose $G$ is a $\mathbb{P}$-generic filter over $\rm{V}$. 
It suffices to verify that for each $\lambda\in K_\mathbb{P}$ greater than  $\delta$, there is $\theta\in K_\mathbb{P}$ and an elementary embedding
$$j\colon V_\lambda[G_\lambda]\rightarrow V_\theta[G_\theta]$$
with $\crit(j)=\delta$ and $j(\delta)>\lambda$.  For since $\lambda $ and $\theta$ are $\mathbb{P}$-reflecting  we then have   $V_{\lambda}[G_{\lambda}]=V[G]_{\lambda}$ and $V_{\theta}[G_{\theta}]=V[G]_{\theta}$. 

So let $\lambda$ be as above. Since $\delta$ is extendible, hence $C^{(1)+}$-extendible, given $\mu\in C^{(1)}$, $\mu>\lambda$, there is $\eta\in C^{(1)}$ and an elementary embedding $j:V_\mu\rightarrow V_\eta$ with $\crit(j)=\delta$ and $j(\delta)>\mu$. Since $K_\mathbb{P}$ is $\Delta_2$-definable, $V_\mu \models \lambda\in K_\mathbb{P}$, hence by elementarity and since $\eta\in C^{(1)}$, $j(\lambda)\in K_\mathbb{P}$.

Let $\theta:=j(\lambda)$. We will show how to lift   $j\restriction V_\lambda :V_\lambda\rightarrow V_\theta$ to an elementary embedding $V_\lambda[G_\lambda]\rightarrow V_\theta[G_\theta]$. Since $\mathbb{P}$ is $\Delta_2$-definable,  $\mathbb{P}^{V_\lambda}=\mathbb{P}\cap V_\lambda=\mathbb{P}_\lambda$, where the rightmost equality follows from $\lambda\in K_\mathbb{P}$. Analogously, the same  holds for $\mathbb{P}^{V_\theta}$. Thus, we have:
\begin{itemize}
\item[(i)] $\mathbb{P}_\lambda\cong \mathbb{P}_\delta\ast \dot{\mathbb{Q}}$, with $|\dot{\mathbb{Q}}|=\lambda$.
\item[(ii)] $\mathbb{P}_\theta\cong \mathbb{P}_{i(\delta)}\ast \dot{\mathbb{Q}}^*$ and $\forces_{\mathbb{P}_{i(\delta)}}\text{$``\dot{\mathbb{Q}}^*$ is weakly homogeneous''}$.
\end{itemize}

\begin{claim}\label{claim2fittingext}
 $\forces_{\mathbb{P}_{i(\delta)}}\text{$``\dot{\mathbb{Q}}^*$ is $i(\delta)$-directed closed''}$.
\end{claim}
\begin{proof}[Proof of claim]
We first show that $\delta$ is an accumulation point of $C(\mathbb{P})$, hence it belongs to $C(\mathbb{P})$. So, fix an ordinal $\sigma_0<\delta$. Since $\mathbb{P}$ is suitable the $\Sigma_3$-formula $\varphi(\sigma_0)$
$$\exists \sigma \forall\rho\forall\alpha\forall X(\sigma_0<\sigma\leq \rho<\alpha\,\wedge \,\alpha\in \overline{K_\mathbb{P}^{(1)}}\,\wedge$$ 
$$X=V_\alpha\rightarrow X\models \text{ $\forces_{\mathbb{P}_\rho}$\text{$``\dot{\mathbb{Q}}_\rho$ is $\sigma_0$-directed closed''})}$$
holds. Since $\delta\in C^{(3)}$ and $\sigma_0<\delta$, there is a witness $\sigma_0<\sigma_1<\delta$ for $\varphi(\sigma_0)$. Arguing inductively we define a sequence $\langle\sigma_n\mid n<\omega\rangle$ of ordinals ${<}\delta$ such that $\sigma_{n+1}$ is a witness for $\varphi(\sigma_n)$. Setting $\sigma_\omega:=\sup_{n<\omega}\sigma_n$ one can easily see that $\sigma_\omega\in C(\mathbb{P})\cap \delta$. 

Since $\delta\in C(\mathbb{P})$, for each $\sigma\in [\delta,\lambda)$, $\forces_{\mathbb{P}_\delta}\text{$``\dot{\mathbb{P}}_{[\delta,\sigma]}$ is $\delta$-directed closed''}$\footnote{Here $\dot{\mathbb{P}}_{[\delta,\sigma]}$ is a $\mathbb{P}_\delta$-name for the iteration $\mathbb{Q}$ up to $\sigma$.}. Since $\lambda$ is inaccessible, $^{<\delta}{V_\lambda}\s V_\lambda$, and therefore the same is true in $V_\lambda$, namely $V_\lambda\models \text{$``\forces_{\mathbb{P}_\delta}$\,\text{$\dot{\mathbb{Q}}$ is $\delta$-directed closed}''}$. By elementarity, $$V_\theta \models \text{$``\forces_{\mathbb{P}_{i(\delta)}}$\,\text{$\dot{\mathbb{Q}}^*$ is $i(\delta)$-directed closed}''}.$$
Since $\theta$ is inaccessible,  $^{<i(\delta)}{V_\theta}\s V_\theta$,  and therefore the quoted sentence displayed above also holds in $V$, as desired.
\end{proof}
From this point on the argument for the lifting of $j\restriction V_\lambda$ to an elementary embedding $V_\lambda[G_\lambda]\rightarrow V_\theta[G_\theta]$ is essentially the same as in Theorem \ref{preservationCn}.
\end{proof}



\begin{theo}\label{keytheonewfitting}
Suppose that $\mathbb{P}$ is a fitting iteration such that $C(\mathbb{P})$ contains all supercompact cardinals, in case there are any. Then forcing with $\mathbb{P}$ preserves \Cn-extendible cardinals, for all $n<\omega$.\footnote{By convention, a cardinal is $C^{(0)}$-extendible iff it is supercompact (cf. page \pageref{conventionsupercompact}).}
\end{theo}
\begin{proof}
We prove the theorem by induction over $n$.  The  case  $n=1$ is covered by   Lemma \ref{2fittingextendibles}.
So, suppose  that $n\geq 2$ and for each $0\leq k< n$ forcing with $\mathbb{P}$ preserves $C^{(k)}$-extendible cardinals. 

We argue similarly as in  the proof of Lemma \ref{2fittingextendibles}. So, let $\delta$ be a \Cn-extendible cardinal and suppose $G$ is a $\mathbb{P}$-generic filter over $\rm{V}$.
Let $\lambda\in K^{(1)}_\mathbb{P}$ be greater than $\delta$, and let $\mu\in C^{(n)}$ be greater than $\lambda$. Then, there is $\eta\in C^{(n)}$ and an elementary embedding
$$j:V_\mu\rightarrow V_\eta$$
with $\crit(j)=\delta$, $j(\delta)>\mu$ and $j(\delta)\in C^{(n)}$. Actually, $j(\delta)$ is $C^{(n-2)}$-extendible. Since $n\geq 2$ and $K_\mathbb{P}$ is $\Delta_{2}$-definable, $V_\mu \models \lambda\in K_\mathbb{P}$. Hence, by elementarity and since $\eta\in C^{(n)}$, $j(\lambda)\in K_\mathbb{P}$. 

Arguing exactly as in the proof of Theorem \ref{2fittingextendibles}, we can lift $j\upharpoonright V_\lambda$ to an elementary embedding
$V[G]_\lambda\rightarrow V[G]_\theta$, where $\theta:=j(\lambda)$.

It only remains to show that $j(\delta)\in C^{(n)}$. Suppose first $n=2$. 
Since $j(\delta)$ is supercompact and every supercompact belongs to $C^{(2)}$, it is enough to show that $j(\delta)$ is supercompact in $V[G]$. By our assumption, $j(\delta)\in C(\mathbb{P})$, hence  $j(\delta)$ is supercompact in $V[G]$ (cf. Theorem \ref{2fittingsupercompact}).

For $1\leq n-2$ and $j(\delta)$ being $C^{(n-2)}$-extendible, our induction hypothesis implies that $i(\delta)$ is $C^{(n-2)}$-extendible in $V[G]$. In particular, $V[G]_{j(\delta)}\prec_{\Sigma_n} V[G]$. Thus,  $\delta$ is $C^{(n)}$-extendible in $V[G]$.
\end{proof}
In the next section we provide several applications of this theorem.

\section{Some applications}

\subsection{Forcing the GCH and related combinatorial principles}

Let $\mathbb{P}=\langle \mathbb{P}_\alpha ; \dot{\mathbb{Q}}_\alpha:\alpha\in \rm{ORD}\rangle$ be the standard Jensen's proper class iteration for forcing the global GCH. Namely, $\mathbb{P}$ is the direct limit of the  iteration with Easton support where $\mathbb{P}_0$ is the trivial forcing and for each ordinal $\alpha$, if $\Vdash_{\mathbb{P}_\alpha}\text{``$\alpha$ is an uncountable cardinal''}$, then $\Vdash_{\mathbb{P}_\alpha}\text{``$\dot{\mathbb{Q}}_\alpha=\mathrm{Add}(\alpha^+,1)$''}$, and $\Vdash_{\mathbb{P}_\alpha}\text{``$\dot{\mathbb{Q}}_\alpha$ is trivial''}$ otherwise. 

\begin{lemma}\label{examplenewfitting}
Assume $\rm{ORD}$ is $\mathbf{\Sigma_2}$-$2$-Mahlo. Then $\mathbb{P}$ is fitting and $C(\mathbb{P})$ contains all supercompact cardinals, in case there are any.
\end{lemma}
\begin{proof}
Clearly, $\mathbb{P}$ is a suitable and weakly homogeneous  iteration.
Also, $\mathbb{P}$ is $\Delta_2$-definable, as ``$p\in \mathbb{P}$'' if and only if $V_\alpha\models ``p\in \mathbb{P}"$, for (some) every   $\alpha \in C^{(1)}$ such that $p\in V_\alpha$. 

Let $K_\mathbb{P}$ be denote the class of all Mahlo cardinals. Clearly, $K_\mathbb{P}$ is $\Delta_2$-definable. Also, since $\mathrm{ORD}$ is $\mathbf{\Sigma}_2$-$2$-Mahlo, $K_\mathbb{P}$ is a $\mathbf{\Sigma_2}$-stationary proper class. 
Let $\kappa\in K_\mathbb{P}$. Thus, $\kappa$ is regular and, clearly, $\mathbb{P}\cap V_\kappa=\mathbb{P}_\kappa$.
\begin{claim}\label{ClaimTsap}
$\kappa$ is $\mathbb{P}$-reflecting.
\end{claim}
\begin{proof}[Proof of claim]
Clearly, $\kappa$ is inaccessible and $\forces_{\mathbb{P}_\kappa} \text{$``\dot{\mathbb{Q}}$ is $\kappa$-distributive''}$, where $\mathbb{P}\cong \mathbb{P}_\kappa\ast \dot{\mathbb{Q}}$. Also, $\mathbb{P}_\kappa\s V_\kappa$. Let $\lambda<\kappa$ be forced by $\mathbb{P}_\kappa$ to be a cardinal. For every inaccessible cardinal $\theta<\kappa$, $\theta\in C(\mathbb{P})$, hence  $$\forces_{\mathbb{P}_\theta}\text{$``\dot{\mathbb{Q}}_{[\theta,\kappa)}$ is $\theta$-distributive''}$$ 
and $\mathbb{P}_\theta=\mathbb{P}\cap V_\theta$. Since $\kappa$ is Mahlo and $\{\theta<\kappa\mid \mathbb{P}_\theta=\mathbb{P}\cap V_\theta\}$ is unbounded, $\mathbb{P}_\kappa$ is $\kappa$-cc. Also, standard arguments show that $\mathbb{P}_\kappa$ forces ``$|\dot{\mathcal{P}}(\lambda)|<\check{\kappa}$''. Altogether, $\mathbb{P}$ forces $\kappa$ to be inaccessible. Finally, we appeal to Proposition \ref{firstprop} to get that $\kappa$ is $\mathbb{P}$-reflecting, as wanted.
\end{proof}
The claim about $C(\mathbb{P})$ containing all supercompacts is obvious.
\end{proof}

In \cite[\S5]{Tsan}, Tsaprounis 
 shows that $\mathbb{P}$ preserves $C^{(n)}$-extendible cardinals. The following theorem  gives  an improvement of his result:

\begin{theo} \hfill
\label{theoTsa}
\begin{enumerate}
\item If $\mathrm{ORD}$ is $\mathbf{\Sigma_2}$-$2$-Mahlo, and $\delta$ is supercompact then
$$\forces_{\mathbb{P}}\text{$``\delta$ is supercompact''}.\footnote{An adhoc argument for  Jensen's iteration can be composed to show that it preserves supercompact cardinals without the assumption, used in Lemma \ref{examplenewfitting}, of $\mathrm{ORD}$ being $\mathbf{\Sigma_2}$-$2$-Mahlo. Nevertheless, the arguments given in the Lemma, also apply to other class forcing iterations considered in the subsections below.}$$
\item If for some $n\geq 1$, $\delta$ is \Cn-extendible then
$$\forces_{\mathbb{P}}\text{$``\delta$ is \Cn-extendible''}.$$
\end{enumerate}
In particular,  if there is a proper class of supercompact cardinals (equivalently, if $\mathrm{VP}(\mathbf{\Pi_1})$ holds),  then forcing with $\mathbb{P}$ preserves \Cn-extendible cardinals, for all $n<\omega$.
\end{theo}
\begin{proof}
This follows from Theorem \ref{keytheonewfitting} and Lemma \ref{examplenewfitting}. The  particular case can be proved using Lemma \ref{OrdMahlo}.
\end{proof}

Recall that a class function $E$ from the class $\mathrm{Reg}$ of infinite regular cardinals to the class of cardinals is an \emph{Easton function} if it satisfies K\"{o}nig's theorem (i.e., $\mathrm{cf}(E(\kappa))>\kappa$, for all $\kappa \in \rm{Reg}$) and is increasingly monotone. Let  $\mathbb{P}_E$ be the direct limit of the iteration $\langle\mathbb{P}_\alpha, \dot{\mathbb{Q}}_\alpha:\,\alpha\in \mathrm{ORD}\rangle$  with Easton support  where $\mathbb{P}_0$ is the trivial forcing and for each ordinal $\alpha$, if $\Vdash_{\mathbb{P}_\alpha}\text{``$\alpha$ is a regular cardinal''}$, then $\Vdash_{\mathbb{P}_\alpha}\text{``$\dot{\mathbb{Q}}_\alpha=\mathrm{Add}(\alpha,E(\alpha))$''}$, and $\Vdash_{\mathbb{P}_\alpha}\text{``$\dot{\mathbb{Q}}_\alpha$ is trivial''}$ otherwise. Standard arguments (\cite[\S15]{Jech}) show that if the GCH holds in the ground model, then $\mathbb{P}_E$ preserves all cardinals and cofinalities and forces  that $2^\kappa=E(\kappa)$, for each $\kappa\in \mathrm{Reg}$. 

\begin{lemma}\label{LemmaPE}
Let $E$ be a $\Delta_2$-definable 
Easton function and assume that $\rm{ORD}$ is $\mathbf{\Sigma_2}$-$2$-Mahlo. 
Then $\mathbb{P}_E$ is fitting and $C(\mathbb{P})$ contains all supercompact cardinals, in case there are any.
\end{lemma}
\begin{proof}
Arguing as in Lemma \ref{examplenewfitting}, $\mathbb{P}_E$ is a $\Delta_2$-definable, weakly homogeneous suitable iteration. Set $\mathbb{P}:=\mathbb{P}_E$ and let $$K_\mathbb{P}:=\{\kappa\mid \text{$\kappa$ is Mahlo and $E[\kappa]\s \kappa$}\}.$$
Clearly, every $\kappa\in K_\mathbb{P}$ is regular and witnesses $\mathbb{P}\cap V_\kappa=\mathbb{P}_\kappa$. The argument for the verification that each $\kappa\in K_\mathbb{P}$ is $\mathbb{P}$-reflecting  is analogous to the   one given  in  Lemma \ref{examplenewfitting}. The claim about $C(\mathbb{P})$ is obvious. Also, since $E$ is $\Delta_2$-definable,  so is $K_\mathbb{P}$. Moreover, $K_\mathbb{P}$ is $\mathbf{\Sigma_2}$-stationary, for if $\mathcal{C}$ is a  $\mathbf{\Sigma_2}$-definable club class of ordinals, then $\mathcal{C}_E:=\{\alpha\in\mathcal{C}\mid E[\alpha]\s \alpha\}$ is a $\mathbf{\Sigma_2}$-definable proper club class, and   since $\rm{ORD}$ is $\mathbf{\Sigma_2}$-$2$-Mahlo, $K_\mathbb{P}\cap \mathcal{C}\neq \emptyset$. 
\end{proof}
Similarly as in Theorem \ref{theoTsa},  we now obtain the following:
\begin{theo}
\label{preserveCnext}
Assume the $\rm{GCH}$ holds and let   $E$ be a $\Delta_2$-definable 
Easton function.\footnote{The GCH assumption is  superfluous  if we are only interested in preserving \Cn-extendible cardinals, $n\geq 1$, for in that case we may first   force   the GCH while preserving \Cn-extendible cardinals (cf. Theorem \ref{theoTsa}).} Then the following hold:
\begin{enumerate}
\item If $\rm{ORD}$ is $\mathbf{\Sigma_2}$-$2$-Mahlo 
and $\delta$ is supercompact then
$$\forces_{\mathbb{P}_E}\text{$``\delta$ is supercompact''}$$
\item If $n\geq 1$ and $\delta$ is \Cn-extendible 
then
$$\forces_{\mathbb{P}_E}\text{$``\delta$ is \Cn-extendible''}$$
\end{enumerate}
In particular, if the $\rm{GCH}$ holds and there is a proper class of supercompact cardinals (equivalently,  $\mathrm{VP}(\mathbf{\Pi_1})$ holds), then
forcing with $\mathbb{P}_E$ preserves \Cn-extendible cardinals, for all $n<\omega$.
\end{theo}

\begin{cor}
Let $E$ be a $\Delta_2$-definable Easton function and   assume the $\rm{GCH}$ holds.  Then the following are true: 
\begin{enumerate}
\item If $\mathrm{VP}(\mathbf{\Pi_{1}})$ holds, then   $\forces_{\mathbb{P}_E}``\mathrm{VP}(\mathbf{\Pi_{1}})$''. 
\item If $\mathrm{VP}({\Pi_{n+1}})$ holds, then   $\forces_{\mathbb{P}_E}``\mathrm{VP}({\Pi_{n+1}})$''. Also, if $\mathrm{VP}(\mathbf{\Pi_{n+1}})$ holds then  $\forces_{\mathbb{P}_E}``\mathrm{VP}(\mathbf{\Pi_{n+1}})$''.
\end{enumerate}
\end{cor}

The next  corollary shows that $\rm{VP}$ is also consistent with any possible   behaviour of the power-set function given by an arbitrary definable Easton function.



\begin{cor}\label{PEandVopenka}
If $\rm{VP}$ holds, then in some class forcing extension, 
for every definable Easton function $E$ there is a further class forcing extension  that preserves $\rm{VP}$ and where $2^\kappa =E(\kappa)$, for every $\kappa\in \mathrm{Reg}$.
\end{cor}

\begin{proof}
First force with the standard Jensen's iteration for forcing the GCH and call the resulting generic extension $V$. By  Corollary \ref{VopenkaIndes}, $V\models \mathrm{VP}$. Then, given a $V$-definable Easton function $E$, force over $V$ with $\mathbb{P}_E$. Once again by Corollary \ref{VopenkaIndes}, $V^{\mathbb{P}_E}\models \rm{VP}$. Finally, since $V\models \mathrm{GCH}$, it follows that $\mathbb{P}_E$ forces $2^\kappa =E(\kappa)$,  for every $\kappa\in\mathrm{Reg}$.  
\end{proof}

\subsection{A remark on Woodin's HOD Conjecture}

The  HOD Dichotomy theorem of Woodin states that if there exists  an extendible cardinal, then either $\rm{V}$ is close to HOD or  is  far from it. Specifically, if $\kappa$ is an extendible cardinal, then either (1): for every singular cardinal $\lambda>\delta$, $\lambda$ is singular in HOD and $(\lambda^+)^{\rm{HOD}}=\lambda^+$, or (2): every regular cardinal $\lambda >\kappa$ is  $\omega$-strongly measurable  in HOD (see  \cite{WOEM1}). Woodin's HOD Hypothesis asserts that there is a proper class of regular cardinals that are not $\omega$-strongly measurable in HOD, and therefore that the first option of the HOD Dichotomy is the true one. Woodin's HOD Conjecture asserts that the HOD Hypothesis is provable in the theory ZFC + ``There exists an extendible cardinal". Our arguments may be used to show that if the HOD Conjecture holds, and therefore it is provable in ZFC + ``There exists an extendible cardinal" that  above the first extendible cardinal every singular cardinal $\lambda$ is singular in HOD and $(\lambda^+)^{\rm{HOD}}=\lambda^+$, there may still be no agreement at all between V and HOD about successors of regular cardinals.  Moreover, many singular cardinals in $\rm{HOD}$ need not be cardinals in $\rm{V}$. Let us give some examples.

\medskip

Let $\mathbb{P}$ be the direct limit of the iteration $\langle \mathbb{P}_\alpha; \dot{\mathbb{Q}}_\alpha:\,\alpha\in \mathrm{ORD}\rangle$ with Easton support,  where $\mathbb{P}_0$ is the trivial forcing and for each ordinal $\alpha$, if $\Vdash_{\mathbb{P}_\alpha}\text{``$\alpha$ is regular''}$ then $\Vdash_{\mathbb{P}_\alpha}\textbf{``$\dot{\mathbb{Q}}_\alpha=\dot{\Coll}(\alpha,\alpha^+)$''}$, and $\Vdash_{\mathbb{P}_\alpha}\text{``$\dot{\mathbb{Q}}_\alpha$ is trivial''}$ otherwise.
Arguing essentially as in Lemma \ref{examplenewfitting} we obtain  the following:

\begin{lemma}\label{woodinisfitting}
Assume $\rm{ORD}$ is $\mathbf{\Sigma}_2$-$2$-Mahlo. Then, $\mathbb{P}$ is  fitting and $C(\mathbb{P})$ contains all supercompact cardinals, in case there are any.
\end{lemma}

\begin{theo}
\label{theoHOD}\hfill
\begin{enumerate}
\item If $\rm{ORD}$ is $\mathbf{\Sigma}_2$-$2$-Mahlo, and $\delta$ is supercompact then
$$\forces_\mathbb{P} \text{$``\delta$ is supercompact''}.$$

\item If $n\geq 1$ and $\delta$ is \Cn-extendible, then
$$\forces_\mathbb{P} \text{$``\delta$ is \Cn-extendible''}$$
\end{enumerate}
Moreover, $\mathbb{P}$ forces $``\forall\lambda\in \mathrm{Reg}\,((\lambda^+)^{\rm{HOD}}<\lambda^+)$''.
\end{theo}

\begin{proof}
For the preservation of $C^{(n)}$-extendible cardinals we combine Lemma \ref{woodinisfitting} and Theorem \ref{keytheonewfitting}. To prove the claim about  successors of regular  cardinals, note that if $\lambda$ is a regular cardinal in $V^{\mathbb{P}}$, then it was also a regular cardinal at stage $\lambda$ of the iteration, hence its successor was collapsed at stage $\lambda +1$. Thus, on the one hand,
$(\lambda^+)^V<(\lambda^+)^{V^{\mathbb{P}}}.$ On the other hand,   $\mathbb{P}$ is weakly homogeneous and ordinal definable,  hence $\rm{HOD}^{{V}^\mathbb{P}}\subseteq HOD^V$ (see, e.g., \cite{Jech} for details). Hence,  in $V^\mathbb{P}$,
$(\lambda^+)^{\rm{HOD}}<\lambda^+$, as wanted.
\end{proof}

\begin{cor}
Forcing with $\mathbb{P}$ preserves $\rm{VP}$ and forces $(\lambda^+)^{\rm{HOD}}<\lambda^+$ for every regular cardinal $\lambda$.
\end{cor}

Theorem \ref{theoHOD} yields the parallel of the main Theorem from \cite{DF}, at the level of $C^{(n)}$-extendible cardinals. 

\medskip

Suppose now that $K$ is a function on the class of  infinite cardinals such that $K(\lambda)> \lambda$, and $K$ is increasingly monotone, for every $\lambda$. 
Let $\mathbb{P}_K$ be the direct limit of an iteration $\langle \mathbb{P}_\alpha; \dot{\mathbb{Q}}_\alpha:\,\alpha\in \mathrm{ORD}\rangle$ with Easton support,  where $\mathbb{P}_0$ is the trivial forcing and for each ordinal $\alpha$, if $\Vdash_{\mathbb{P}_\alpha}\text{``$\alpha$ is regular''}$ then $\Vdash_{\mathbb{P}_\alpha}\textbf{``$\dot{\mathbb{Q}}_\alpha=\dot{\Coll}(\alpha,K(\alpha))$''}$, and $\Vdash_{\mathbb{P}_\alpha}\text{``$\dot{\mathbb{Q}}_\alpha$ is trivial''}$ otherwise. 
Notice that for each $\alpha$ such that $\Vdash_{\mathbb{P}_\alpha}\text{``$\alpha$ is regular''}$, the remaining part of the iteration after stage $\alpha$ is $\alpha$-closed, hence it preserves $\alpha$. Also note that if $K$ is $\Delta_m$-definable ($m\geq 1$), then $\mathbb{P}_K$ is also $\Delta_m$-definable.  Clearly, $\mathbb{P}_K$ is suitable and weakly homogeneous. 

\begin{lemma}\label{woodinisfitting2}
Assume $\rm{ORD}$ is $\mathbf{\Sigma}_2$-$2$-Mahlo. Let $K$ be a $\Delta_2$-definable class function as above. Then, $\mathbb{P}_K$ is a  fitting iteration and $C(\mathbb{P})$ contains all supercompact cardinals, in case there are any.
\end{lemma}
\begin{proof}
Let $K_\mathbb{P}:=\{\kappa\mid \text{$\kappa$ is Mahlo and $K[\kappa]\s \kappa$}\}$. The proof that $K_\mathbb{P}$ is a witness for the fittingness of $\mathbb{P}$ is similar to those given in lemmas \ref{LemmaPE} and \ref{woodinisfitting}. The claim about $C(\mathbb{P})$ is obvious.
\end{proof}

\begin{theo}
\label{preserveCnext2}
Let $K$ be a  $\Delta_2$-definable class function as above. 

\begin{enumerate}
\item If $\rm{ORD}$ is $\mathbf{\Sigma_2}$-$2$-Mahlo, and $\delta$ is supercompact then
$$\forces_{\mathbb{P}_K} \text{$``\delta$ is supercompact.''}$$
\item  If $\delta$ is \Cn-extendible, for some $n\geq 1$, then
$$\forces_{\mathbb{P}_K} \text{$``\delta$ is \Cn-extendible.''}$$
\end{enumerate}
 Moreover, $\mathbb{P}_K$ forces $$(\lambda^+)^{\rm{HOD}}\leq K(\lambda)<\lambda^+$$ for all infinite regular cardinals $\lambda$.
\end{theo}

\begin{proof}
The preservation of \Cn-extendible cardinals, $n<\omega$, follows from Lemma \ref{woodinisfitting2} and Theorem \ref{keytheonewfitting}.

If $G$ is $\mathbb{P}_K$-generic over $\rm{V}$ and $\lambda$ is regular in $V[G]$, then it is also regular at the $\lambda$-stage  of the iteration. Hence, $\mathbb{Q}_\lambda =\Coll(\lambda ,K(\lambda))$, and therefore $K(\lambda)<\lambda^+$ holds in $V[G]$. The other inequality (i.e., $(\lambda^+)^{\rm{HOD}}\leq K(\lambda)$) follows from the fact that $\mathbb{P}_K$ is weakly homogeneous and ordinal definable, and thus that $\rm{HOD}^{V[G]}\subseteq \rm{HOD}^V$.
\end{proof}

The theorem above implies that many kinds of disagreement between successors of regulars in $\rm{HOD}$ and in $\rm{V}$ may be forced while preserving $C^{(n)}$-extendible cardinals. It also implies that one can destroy many singular cardinals in $\rm{HOD}$ while preserving $C^{(n)}$-extendible cardinals. For example, let  $K$ be such that $K(\lambda)$ is the least singular cardinal in $\rm{HOD}$ greater than $\lambda$, i.e., $K(\lambda)=(\lambda^{+\omega})^{\rm{HOD}}$. It is easily seen that $K$, and therefore also $\mathbb{P}_K$ as defined above, are $\Delta_2$-definable. Then we have the following.  

\begin{cor}
For each $n\geq 1$, $\mathbb{P}_K$ preserves $C^{(n)}$-extendible cardinals  and forces 
$``\forall\lambda\in\mathrm{Reg}\,((\lambda^{+\omega})^{\rm{HOD}}<\lambda^+)"$.
\end{cor}

\subsection{On diamonds}

Other combinatorial statements that we can force while preserving \Cn-extendible cardinals  are  the diamond principles $\diamondsuit_S$.
Namely, given a stationary set $S\s \kappa$, a sequence $\langle A_\alpha:\,\alpha\in S\rangle$ is a $\diamondsuit_S$\textit{-sequence} if $A_\alpha\subseteq \alpha$ and for every $A\subseteq\kappa$ the set $\{\alpha\in S:\, A\cap\alpha=A_\alpha\}$ is stationary. We say that $\diamondsuit_S$ holds if there is a $\diamondsuit_S$\textit{-sequence}. 

\medskip

It is well-known that $\mathrm{Add}(\kappa^+,1)$ automatically forces $\diamondsuit_{S}$, for every stationary  $S\subseteq\kappa^+$ in $V^{\mathrm{Add}(\kappa^+,1)}$. Thus, from Theorem \ref{theoTsa}, we obtain:
 
 \begin{cor}
\hfill

\begin{enumerate}
\item If $\rm{ORD}$ is $\mathbf{\Sigma_2}$-$2$-Mahlo and $\delta$ is a  supercompact cardinal, then there is a generic extension where  $\delta$ is still supercompact and $\diamondsuit_{S}$ holds for every cardinal $\kappa$ and every stationary  $S\subseteq\kappa^+$.
\item If $n\geq 1$ and $\delta$ is a  \Cn-extendible cardinal,  there is a generic extension where  $\delta$ is \Cn-extendible and $\diamondsuit_{S}$ holds for every cardinal $\kappa$ and every stationary  $S\subseteq\kappa^+$.
\end{enumerate}
Hence, if $\rm{VP}$ holds in $\rm{V}$,     there is a generic extension where $\rm{VP}$ holds together with $\diamondsuit_{S}$, for every $\kappa$ and every stationary  $S\subseteq\kappa^+$.
 \end{cor}
 
Another relevant diamond principle   is the so called $\diamondsuit^+_{\kappa^+}$-principle.

A sequence $\langle \mathcal{A}_\alpha:\,\alpha\in\kappa^+\rangle$ is  a $\diamondsuit^+_{\kappa^+}$-sequence if $\mathcal{A}_\alpha\in [\mathcal{P}(\alpha)]^{\leq \kappa}$ and for every $A\subseteq\kappa^+$ there is a club $C\subseteq\kappa^+$ such that
 $$C\subseteq\{\alpha\in \kappa^+\mid \, A\cap \alpha\in \mathcal{A}_\alpha\,\wedge\, C\cap \alpha\in \mathcal{A}_\alpha\}.$$
We say that $\diamondsuit^+_{\kappa^+}$ holds if there is a $\diamondsuit^+_{\kappa^+}$-sequence.  

In \cite[Theorem 12.2]{CummingsSquare} it is shown  that,   assuming $2^\kappa=\kappa^+$ and $2^{\kappa^+}=\kappa^{++}$, there is a $\kappa^+$-closed and $\kappa^{++}$-cc forcing notion that forces $\diamondsuit^+_{\kappa^+}$. The forcing is  an iteration $\mathbb{D}^+_{\kappa^{++}}=\langle \mathbb{P}_\alpha,\dot{\mathbb{Q}}_\beta: \beta<\alpha\leq\kappa^{++}\rangle$ with supports of size $\leq\kappa$, where $\mathbb{P}_0$ is the natural forcing notion that introduces a sequence $\vec {\mathcal{A}}$ of the right form 
whereas the rest of the iterates forces the  club sets $C\subseteq\kappa^+$ witnessing that $\vec{\mathcal{A}}$ is  a $\diamondsuit^+_{\kappa^+}$-sequence. 

\smallskip

Arguing as in Lemma~\ref{examplenewfitting} one has the following: 

\begin{lemma}\label{Disfitting}
Assume  the $\rm{GCH}$ holds and that $\rm{ORD}$ is $\mathbf{\Sigma_2}$-$2$-Mahlo. Let $\mathbb{D}$ be the standard Easton support iteration   of the forcings $\mathbb{D}^{+}_{\kappa^{++}}$, for $\kappa$ a cardinal. Then $\mathbb{D}$ is fitting and $C(\mathbb{P})$ contains all supercompact cardinals, in case there are any.
\end{lemma}





\begin{theo}
Assume the $\rm{GCH}$ holds.
\begin{enumerate}
\item If $\rm{ORD}$ is $\mathbf{\Sigma_2}$-$2$-Mahlo and $\delta$ is a  supercompact cardinal, then forcing with $\mathbb{D}$ preserves the supercompactness of    $\delta$ and yields a generic extension where $\diamondsuit^+_{\kappa^+}$ holds, for every cardinal $\kappa$.
\item If $n\geq 1$ and $\delta$ is a  \Cn-extendible cardinal, then forcing with $\mathbb{D}$ preserves the  \Cn-extendibility of $\delta$ and yields a generic extension where $\diamondsuit^+_{\kappa^+}$ holds, for every cardinal $\kappa$.
\end{enumerate}
Hence,  
if $\rm{VP}$ and the $\rm{GCH}$ hold in $\rm{V}$, then $\rm{VP}$ also holds in $V^{\mathbb{D}}$, together with $\diamondsuit_{\kappa^+}^+$, for every cardinal $\kappa$.  
\end{theo}
The claim of the theorem above referring to $\rm{VP}$ was first  proved by Brooke-Taylor in \cite[Corollary 26]{Broo}


\subsection{On weak square sequences}
A classical result due to Solovay is that Jensen's square principle $\square_\lambda$ must fail for every cardinal $\lambda$ greater than or equal to the first strongly compact cardinal \cite{SolGCH}. This result  was subsequently refined by Jensen, who proved that $\square_\lambda$ fails for any subcompact cardinal $\lambda$, a much weaker notion than supercompactness. Further sharper results are due to Brooke-Taylor and Friedman \cite{BrooFrie} and to Bagaria and Magidor \cite{BM2}. In our context, namely with the existence of extendible cardinals, and therefore with the failure of $\square_\lambda$ for a tail of $\lambda$'s, we shall consider $\square_{\lambda,\mu}$-principles, a weak form of the square principle introduced by Schimmerling in \cite{Schi}.

\smallskip

Following up on Solovay's work, Shelah proved that if $\kappa$ is supercompact and $\rm{cof}(\lambda)<\kappa <\lambda$, then $\square_{\lambda,\lambda}$ (also known as  $\square^*_\lambda$) fails \cite[\S2]{CummingsSquare}. Also, Burke and Kanamori showed  that if $\kappa$ is $\lambda^+$-strongly compact, then $\square_{\lambda,<\cof(\lambda)}$ fails  \cite{CummingsSquare}.\footnote{Note that this is interesting only for $\cof(\lambda)\geq \kappa$.} The remaining cases, namely $\square_{\lambda,\mu}$ with $\kappa\leq \cof(\lambda)\leq\mu\leq\lambda$ turned out to be consistent with the existence of  a supercompact cardinal $\kappa$ \cite[\S9]{CummingsSquare}. 
Specifically, in \cite[Theorem 9.2]{CummingsSquare}, the authors define, for each $\kappa\leq\cof(\lambda)<\lambda$, a forcing notion $\mathbb{S}_\lambda$   which forces $\square_{\lambda, \cof(\lambda)}$ and  is $\cof(\lambda)$-directed closed and $\lambda$-strategically closed. Thus, if $\kappa$ is  Laver-indestructible, then forcing with $\mathbb{S}_\lambda$ preserves the supercompactness of $\kappa$. In this section we shall extend this result to  \Cn-extendible cardinals. 


\medskip

Let $K$ be the class function with $\dom(K)=\mathrm{Card}\setminus \aleph_1$, such that $K(\lambda):= ``\text{The  first singular cardinal of cofinality $\lambda^+$''}$. Observe that $K$ is $\Delta_2$-definable. Now let $\mathbb{P}_K$ be denote the direct limit of $\langle \mathbb{P}_\alpha;\dot{\mathbb{Q}}_\alpha:\,\alpha\in \rm{ORD}\rangle$, the iteration with Easton support where $\mathbb{P}_0$ is the trivial forcing and for every $\alpha$, if $\Vdash_{\mathbb{P}_\alpha}\text{``$\alpha\in\dot{\mathrm{Card}}\setminus \dot{\aleph}_1$''}$, then $\Vdash_{\mathbb{P}_\alpha}\text{``$\dot{\mathbb{Q}}_\alpha=\dot{\mathbb{S}}_{K(\alpha)}$''}$, 
and $\Vdash_{\mathbb{P}_\alpha}\text{``$\dot{\mathbb{Q}}_\alpha=\{\mathds{1}\}$''}$, otherwise. 

\begin{remark}
 For every $\alpha$, if $\Vdash_{\mathbb{P}_\alpha}\text{$``\alpha\in\dot{\mathrm{Card}}\setminus \dot{\aleph}_1$"}$, then $\mathbb{P}_{[\alpha,\rm{ORD})}$ is $\alpha^+$-directed closed and $K(\alpha)$-strategically closed.
\end{remark}

We now prove that $\mathbb{P}_K$ forces class many instances of $\square_{K(\lambda),\lambda^+}$.

\begin{lemma}\label{PKdoesthejob}
  $\forces_{\mathbb{P}}\text{$``\forall \lambda\geq \dot{\aleph}_1\,(\square_{K(\lambda),\lambda^+}$ holds)''}$.
\end{lemma}
\begin{proof}
Suppose $\lambda$ is an uncountable cardinal in $V^\mathbb{P}$. Then, by the above remark, $\lambda$ is also an uncountable cardinal in $V^{\mathbb{P}_\lambda}$, hence \linebreak$\Vdash_{\mathbb{P}_\lambda}``\dot{\mathbb{Q}}_\lambda=\dot{\mathbb{S}}_{K(\lambda)}"$. Thus, $\square_{K(\lambda),\lambda^+}$ holds  in $V^{\mathbb{P}_{\lambda+1}}$. Again, by the above remark, $\mathbb{P}_{[\lambda+1,\rm{ORD})}$ is distributive enough to preserve $\square_{K(\lambda),\lambda^+}$.  
\end{proof}

Arguing as in Lemma~\ref{LemmaPE} we can show that $\mathbb{P}_K$ is a fitting iteration. Precisely, we have the following:

\begin{lemma}\label{squareisfitting}
Assume $\rm{ORD}$ is $\mathbf{\Sigma_2}$-$2$-Mahlo. Then  $\mathbb{P}_K$ is fitting and $C(\mathbb{P})$ contains all supercompact cardinals, in case there are any.
\end{lemma}


\begin{theo}\label{squaretheorem}
\hfill
\begin{enumerate}
\item If $\rm{ORD}$ is $\mathbf{\Sigma_2}$-$2$-Mahlo,  then $\mathbb{P}_K$ preserves supercompact cardinals.
\item For each $n\geq 1$, forcing  with $\mathbb{P}_K$ preserves \Cn-extendible cardinals.
\end{enumerate}
Moreover, in any generic extension by $\mathbb{P}_K$ the following holds: for every uncountable cardinal $\lambda$, if $K(\lambda)$ is the first singular cardinal of cofinality $\lambda^+$, then  $\square_{K(\lambda),\lambda^+}$ holds.
\end{theo}

\begin{cor}\label{VPSquares}
If  $\VP$ holds, then there is a class forcing iteration that preserves $\VP$ and forces  $\square_{\lambda,\cof(\lambda)}$, for a proper class of singular cardinals $\lambda$.
\end{cor}

\section{General class forcing iterations}
In this section we follow up the discussion at the end of Section 6 about non weakly homogeneous suitable iterations. One prominent example  is the iteration $\mathbb{P}$ that forces $\rm{V}=\rm{HOD}$ by coding the universe into the power-set function pattern. This iteration is suitable but not weakly homogenous. One may also want to consider class forcing iterations $\mathbb{P}$ over some model $M$ such that $\mathbb{P}$ is not definable in $M$. To deal  with such general class forcing notions  we shall work within the theory  $\rm{ZFC}_\mathbb{P}$, namely  ZFC with the axiom schemata of Separation and Replacement allowing for formulas in the  language of set theory with the additional predicate symbol $\mathbb{P}$.   


\begin{defi}[$\mathbb{P}$-\,$C^{(n)}$-extendible cardinal]\label{DefiPCnextendible}
For $n\geq 1$, we say that a cardinal $\delta$ is \emph{$\mathbb{P}$-\,$C^{(n)}$-extendible} if for every cardinal  $\lambda\in C^{(n)}_{\mathbb{P}}$, $\lambda>\delta$, there is an ordinal $\theta$ and an elementary embedding $$j: \langle V_\lambda,\in, \mathbb{P}\cap V_\lambda\rangle\rightarrow \langle V_\theta,\in, \mathbb{P}\cap V_\theta\rangle$$ with $\crit(j)=\delta$, $j(\delta)>\lambda$, and $j(\delta)\in C^{(n)}$. 
 If, moreover, we can pick $\theta\in C^{(n)}_{\mathbb{P}}$, then we  say that $\delta$ is $\mathbb{P}$-\,$C^{(n)+}$-extendible.
\end{defi}

Similarly,   we may also consider the notion of $\mathbb{P}$-$\Sigma_n$-supercompactness, for a general  class $\mathbb{P}$ which is  not necessa\-rily definable. 

\begin{defi}[$\mathbb{P}$-$\Sigma_n$-supercompactness]\label{PSigmaSupercom2}
If $n\geq 1$,   then  a cardinal $\delta$ is  $\mathbb{P}$-$\Sigma_n$-supercompact if  for every  $\lambda\in C^{(n)}_\mathbb{P}$ greater than $\delta$,  and every $a\in V_\lambda$ there exist $\bar{\delta}<\bar{\lambda}<\delta$ and $\bar{a}\in V_{\bar{\lambda}}$, and there exists an elementary embedding $j: V_{\bar{\lambda}}\longrightarrow V_\lambda$ such that:
\begin{itemize}
\item $\crit(j)=\bar{\delta}$ and $j(\bar{\delta})=\delta$.
\item $j(\bar{a})=a$.
\item $\bar{\lambda}\in C^{(n)}_\mathbb{P}$.
\end{itemize}
\end{defi}

The same arguments as in the proof of Theorem \ref{MagidorLikecharacteri} yield  the following equivalence.

\begin{theo}
\label{theoP+}
For every $n\geq 1$, every class $\mathbb{P}$, and every cardinal $\kappa$, the following  are equivalent:
\begin{enumerate}
\item $\kappa$ is $\mathbb{P}$-\,$C^{(n)}$-extendible.
\item $\kappa$ is $\mathbb{P}$-\,$\Sigma_{n+1}$-supercompact.
\item $\kappa$ is $\mathbb{P}$-\,$C^{(n)+}$-extendible.\footnote{See Remark~\ref{RemarkOnC+extendibility}.} 
\end{enumerate}
\end{theo}
We will say that a cardinal $\delta$ is \emph{$\mathbb{P}$-$C^{(n)}$-extendible with $\mathbb{P}$-$\Sigma_n$-reflecting target} if in Definition \ref{DefiPCnextendible} we may moreover require that $j(\delta)$ is $\mathbb{P}$-$\Sigma_n$-reflecting (cf. Definition \ref{Psigmanreflecting}). Likewise, one defines the notion of \emph{$\mathbb{P}$-$C^{(n)+}$-extendible with $\mathbb{P}$-$\Sigma_n$-reflecting target}. Arguing as usual one can check that both notions are equivalent.

Clearly, any $\mathbb{P}$-$C^{(n)+}$-extendible with $\mathbb{P}$-$\Sigma_n$-reflecting target is \Cn-extendi\-ble. Moreover, if the predicate $\mathbb{P}$ is definable with low complexity and satisfies some minor requirements, then, as the next proposition shows, a cardinal is $\mathbb{P}$-$C^{(n)+}$-extendible with $\mathbb{P}$-$\Sigma_n$-reflecting target if and only if is $C^{(n)}$-extendible.
\begin{prop}\label{PredicateandCnext}
Let $n,m\geq 1$ with $m\leq n$. Let  $\mathbb{P}$ be a $\Delta_{m+1}$-definable $\ord$-length forcing iteration and assume that  $$\{\kappa\mid \mathbb{P}\cap V_\kappa=\mathbb{P}_\kappa\,\wedge\,\text{$\kappa$ is $\mathbb{P}$-reflecting}\}$$ is a proper class containing  all $C^{(n-2)}$-extendible cardinals.\footnote{By convention, a $C^{(-1)}$-extendible cardinal is a Mahlo cardinal.}  Then, \linebreak every $C^{(n)}$-extendible cardinal belonging to $C^{(m+n+1)}$ is $\mathbb{P}$-$C^{(n)+}$-exten\-dible with $\mathbb{P}$-$\Sigma_n$-reflecting target.

 In particular, if $\mathbb{P}$ is $\Delta_2$-definable, a cardinal is \Cn-extendible if and only if it is   $\mathbb{P}$-$C^{(n)+}$-extendible with $\mathbb{P}$-$\Sigma_n$-reflecting target. 
\end{prop}
\begin{proof}
Let $\delta$ be a $C^{(n)}$-extendible cardinal in $C^{(m+n+1)}$. 
\begin{claim}\label{auxiliaryclaimextendibles}
$\delta$ is $\mathbb{P}$-\Cn-extendible. 
\end{claim}
\begin{proof}[Proof of claim]
Let $\mu>\lambda>\delta$ be with $\lambda\in C^{(n)}_\mathbb{P}$ and $\mu\in C^{(n)}$. Since $\delta$ is actually $C^{(n)+}$-extendible (cf. Remark \ref{RemarkOnC+extendibility}) we may pick $\eta\in C^{(n)}$ together with an elementary embedding $j\colon V_\mu\rightarrow V_{\eta}$ with $\crit(j)=\delta$, $j(\delta)>\mu$ and $j(\delta)\in C^{(n)}$. Observe that  $V_{j(\lambda)}\prec_{\Sigma_n} V_\eta\prec_{\Sigma_n} V$, hence $$j(\langle V_\lambda,\in,\mathbb{P}^{V_\lambda}\rangle)=\langle V_{j(\lambda)},\in, \mathbb{P}^{V_{j(\lambda)}}\rangle=\langle V_{j(\lambda)},\in, \mathbb{P}\cap V_{j(\lambda)}\rangle,$$
where the right-most equality follows from $\Delta_{m+1}$-definability of $\mathbb{P}$ and $m\leq n$.  Altogether,  $j\upharpoonright\langle V_\lambda,\in,\mathbb{P}\cap V_\lambda\rangle$ yields the desired embedding.
\end{proof}
By Theorem \ref{theoP+}, $\delta$ is actually $\mathbb{P}$-$C^{(n)+}$-extendible, so we concentrate on proving the other assertion. Let $\lambda\in C^{(n)}_\mathbb{P}$ with $\lambda>\delta$, together with an elementary embedding
$$j\colon \langle V_\lambda,\in,\mathbb{P}\cap V_\lambda\rangle\rightarrow \langle V_\theta,\in,\mathbb{P}\cap V_\theta\rangle,$$
with $\crit(j)=\delta$, $\theta\in C^{(n)}_\mathbb{P}$, $j(\delta)>\lambda$ and $j(\delta)\in C^{(n)}$. Note that $j(\delta)$ is $C^{(n-2)}$-extendible hence $\mathbb{P}$-reflecting  and a witness for $\mathbb{P}\cap V_{j(\delta)}=\mathbb{P}_{j(\delta)}$. 
\begin{claim}\label{secondauxiliaryclaimcnext}
$j(\delta)\in C^{(n)}_\mathbb{P}$.
\end{claim}
\begin{proof}[Proof of Claim]
Since  $\mathbb{P}$ is $\Delta_{m+1}$-definable, $C^{(n)}_\mathbb{P}$ is a $\Pi_{n+m}$-definable club class (cf. Proposition \ref{complexityclubclass}). In particular, as $\delta\in C^{(n+m+1)}$, $\delta$ is an accumulation point of $C^{(n)}_\mathbb{P}$ and thus a member of $C^{(n)}_\mathbb{P}$.

 Since $\mathbb{P}$ is $\Delta_{m+1}$, $m\leq n$ and $\delta,\lambda\in C^{(n)}_\mathbb{P}$,  $\langle V_\delta,\in, \mathbb{P}^{V_\delta}\rangle\prec_{\Sigma_n}\langle V_\lambda,\in,\mathbb{P}^{V_\lambda}\rangle$. By elementarity, $\langle V_{j(\delta)},\in, \mathbb{P}^{V_{j(\delta)}}\rangle\prec_{\Sigma_n}\langle V_\theta,\in,\mathbb{P}^{V_\theta}\rangle$ and so, since $\theta\in C^{(n)}_\mathbb{P}$, $\langle V_{j(\delta)},\in, \mathbb{P}^{V_{j(\delta)}}\rangle\prec_{\Sigma_n}\langle V,\in,\mathbb{P}\rangle$. Once again, since $\mathbb{P}$ is $\Delta_{m+1}$-definable, $m\leq n$  and $j(\delta)\in C^{(n)}$, $\mathbb{P}^{V_{j(\delta)}}=\mathbb{P}\cap V_{j(\delta)}$, hence  $j(\delta)\in C^{(n)}_\mathbb{P}$.
\end{proof}
The above claims combined give the proof of the proposition.
\end{proof}

%

\medskip

The main result of the section is the following:
\begin{theo}\label{NonHomoVopenka} 
Let $\mathbb{P}$ be a (not necessarily definable) suitable iteration. Assume that there is a proper class of  $\mathbb{P}$-reflecting cardinals $\lambda$ such that $\mathbb{P}_\lambda =\mathbb{P}\cap V_\lambda$.\footnote{In particular, $\mathbb{P}$ is adequate (cf. Definition \ref{defadequate}).} For each $n\geq 1$, if $\delta$ is    $\mathbb{P}$-\,$C^{(n)}$-extendible with $\mathbb{P}$-$\Sigma_n$-reflecting target then $$\forces_\mathbb{P}\text{$``\delta$ is $C^{(n)}$-extendible''}.$$
Actually for $n=1$ the above is true by just assuming  that $\delta$ is $\mathbb{P}$-$C^{(1)}$-extendible.
\end{theo}
\begin{proof} 
Let $\lambda >\delta$ be $\mathbb{P}$-reflecting and such that $\mathbb{P}_\lambda =\mathbb{P}\cap V_\lambda$. 
It will be sufficient to prove that if $G_\lambda$ is $\mathbb{P}_\lambda$-generic over $\rm{V}$, then in the generic extension $V[G_\lambda]$, the set $D$ of conditions $r\in \mathbb{P}_{[\lambda,Ord)}$ that  force the existence of an elementary embedding 
$$j:V[G_\lambda][\dot{G}_{[\lambda, Ord)}]_\lambda\to  V[G_\lambda][\dot{G}_{[\lambda, Ord)}]_\theta$$ some $\theta$, with  $\crit(j)=\delta$, $j(\delta)>\lambda$,  and $j(\delta)\in C^{(n)}$,  
is dense in $\mathbb{P}_{[\lambda,Ord)}$.

In $V[G_\lambda]$, let $r$ be a condition in $\mathbb{P}_{[\lambda,Ord)}$. Back in $\rm{V}$, let  $\mu \in C^{(n)}_{\mathbb{P}}$  be greater than $\lambda$  and such that 
$$ \Vdash_{\mathbb{P}_\mu}\text{``$\mathbb{P}_{[\mu,Ord)}$ is $|\mathbb{P}_\lambda|^+$-directed closed''}.$$

Since $\delta$ is $\mathbb{P}$-\,$C^{(n)+}$-extendible with $\mathbb{P}$-$\Sigma_n$-reflecting target, in the ground model $\rm{V}$ there exists an elementary embedding $$j: \langle V_\mu, \in ,\mathbb{P}\cap V_\mu\rangle\rightarrow  \langle V_\theta, \in ,\mathbb{P}\cap V_\theta\rangle$$ with $\crit(j)=\delta$ such that $j(\delta)>\mu$,  $\theta\in C^{(n)}_\mathbb{P}$ and $ j(\delta)$ being $\mathbb{P}$-$\Sigma_n$-reflecting.

For each $q\in \mathbb{P}_\lambda$ there is an ordinal $\alpha<\delta$ such that $\supp(q)\cap \delta\subseteq\alpha$. Hence, $\supp(j(q))\cap j(\delta)\subseteq \alpha$,  and so $j(q)$ is a $\mathbb{P}_{j(\lambda)}$-condition such that
$$j(q)(\beta)=\left\{\begin{array}{ccc}
q(\beta)& \text{if}& \beta<\alpha.\\
\mathds{1}& \text{if}& \beta\in [\alpha, j(\delta)).
\end{array}\right.$$ 
Since $\mu<j(\delta)$ we have that $\supp(j(q))\cap [\lambda,\mu) =\emptyset$. So, by our choice of the ordinal $\mu$, in $V[G_\lambda]$ we can take $r^*\in\mathbb{P}_{[\mu, Ord)}$ such that 
$$ \Vdash_{\mathbb{P}_{[\lambda ,\mu)}}``r^*\leq j(q)\upharpoonright [\mu, j(\lambda))"$$ for all $q\in G_\lambda$. Then, the condition $r\wedge r^\ast$ such that
$$(r\wedge r^*)(\beta):=\left\{\begin{array}{ccc}
r(\beta)& \text{if}& \beta\in [\lambda,\mu).\\
r^*(\beta)& \text{if}& \beta\in[\mu, j(\lambda)).
\end{array}\right.$$ 
is well-defined and works as a master condition for $j$ and the forcing $\mathbb{P}_{j(\lambda)}/G_\lambda$, because  $$r\wedge r^*\Vdash_{\mathbb{P}_{j(\lambda)}/G_\lambda} j``G_\lambda\subseteq \dot{G}_{j(\lambda)}.$$ 
Thus, for any $\mathbb{P}_{j(\lambda)}$-generic filter $G_{j(\lambda)}$ over $\rm{V}$ extending $G_\lambda$ and containing $r\wedge r^\ast$, the  elementary embedding
$$j\restriction V_\lambda : \langle V_\lambda, \in , \mathbb{P}\cap V_\lambda  \rangle \to \langle V_{j(\lambda)}, \in , \mathbb{P}\cap V_{j(\lambda)}  \rangle$$
lifts to an elementary embedding 
$$j^\ast: \langle V_\lambda[G_\lambda], \in , \mathbb{P}\cap V_\lambda[G_\lambda]  \rangle \to \langle V_{j(\lambda)}[G_{j(\lambda)}], \in , \mathbb{P}\cap V_{j(\lambda)}[G_{j(\lambda)}]  \rangle.$$
Now, since $\lambda$ is $\mathbb{P}$-reflecting, $\mathbb{P}$ forces that $V_\lambda[\dot{G}_\lambda]=V[\dot{G}]_\lambda$. Hence, by the choice of $\mu$, the same is forced by $\mathbb{P}_\mu$.  By the elementarity of $j$, the structure $\langle V_\theta, \in ,\mathbb{P}\cap V_\theta \rangle$ thinks that the forcing $\mathbb{P}\cap V_\theta$ forces  $V_{j(\lambda)}[\dot{G}_{j(\lambda)}]=V[\dot{G}]_{j(\lambda)}$. So, since $\theta \in C^{(n)}_{\mathbb{P}}$, $\mathbb{P}$ forces the same. 
Also, since $j(\delta)$ is $\mathbb{P}$-$\Sigma_n$-reflecting in $\rm{V}$ and $V_{j(\delta)}\models \mathrm{ZFC}$,  Lemma \ref{PsigmareflectingCncardinal} yields $$ \forces_\mathbb{P}``j^*(\delta)\in\dot{C}^{(n)}\text{''.}$$

 We have thus found a condition below $r$, namely $r\wedge r^\ast$, forcing  the existence of an elementary embedding 
$$j^\ast: \langle V[\dot{G}]_\lambda, \in , \mathbb{P}\cap V[\dot{G}] _\lambda\rangle \to \langle V[\dot{G}]_{j(\lambda)}, \in , \mathbb{P}\cap V[\dot{G}]_{j(\lambda)}  \rangle$$
with $\crit(j^\ast)=\delta$, $j^\ast(\delta)>\lambda$, and $j^\ast(\delta)\in C^{(n)}$, as wanted.
\end{proof}

\subsection{$\VP$ and non-homogeneous suitable iterations}\label{SectionVopenkaNoHomo}

We now use  Theorem \ref{NonHomoVopenka} to prove Theorem \ref{AlmostVopenka} without the homogeneity assumption on $\mathbb{P}$, hence yielding the desired refinement of Brooke-Taylor's main theorem of \cite{Broo}.

\begin{theo}\label{VopenkaNonHom}
Let $n,m\geq 1$ with $m\leq n$. Let $\mathbb{P}$ be a  $\Delta_{m+1}$-definable suitable iteration. 
If  $\mathrm{VP}(\mathbf{\Pi_{m+n+1})}$ holds then $V^\mathbb{P}\models\mathrm{VP}(\mathbf{\Pi_{n+1}})$. 

In particular, if   $\mathrm{VP}$ holds and $\mathbb{P}$ is a definable suitable iteration,  then $$V^\mathbb{P}\models \mathrm{VP}.$$
\end{theo}
\begin{proof}
By Theorem \ref{VopenkaBagaria}, $\mathrm{VP}(\mathbf{\Pi_{m+n+1})}$ yields the existence of a proper class of $C^{(m+n)}$-extendible cardinals. Also, by  Theorem \ref{VopParam}, it entails the existence of a proper class of  $\mathbb{P}$-$\Sigma_{n+1}$-reflecting cardinals, hence a proper class of $\mathbb{P}$-reflecting cardinals $\lambda$ such that $\mathbb{P}_\lambda=\mathbb{P}\cap V_\lambda$. In particular, the assumptions of Theorem \ref{VopenkaNonHom} are met.
\begin{claim}
For each $n\geq 2$, every  $C^{(m+n)}$-extendible is $\mathbb{P}$-$C^{(n)}$-extendible with $\mathbb{P}$-$\Sigma_n$-reflecting target.
\end{claim}
\begin{proof}[Proof of claim]
If $\delta$ is $C^{(m+n)}$-extendible then it is \Cn-extendible and $\Sigma_{m+n+2}$-correct. Thus, Claim \ref{secondauxiliaryclaimcnext}  implies that $\delta$ is  $\mathbb{P}$-$C^{(n)}$-extendible.

Let $\lambda\in C^{(n)}_\mathbb{P}$ be with $\lambda>\theta$. By $\mathbb{P}$-$C^{(n)+}$-extendibility of $\delta$, there is $\theta\in C^{(n)}_\mathbb{P}$ and an elementary embedding
$$j: \langle V_\lambda,\in,\mathbb{P}\cap V_\lambda\rangle\rightarrow \langle V_\theta,\in,\mathbb{P}\cap V_\theta\rangle,$$
with $\crit(j)=\delta$, $j(\delta)>\lambda$ and $j(\delta)\in C^{(n)}$. 

\begin{subclaim}
$j(\delta)$ is $\mathbb{P}$-$\Sigma_n$-reflecting. 
\end{subclaim}
\begin{proof}[Proof of subclaim]
Since $\delta$ is $C^{(n+m)}$-extendible, hence \Cn-extendible and $\delta\in C^{(n+m+1)}$,  Claim \ref{secondauxiliaryclaimcnext}  yields $j(\delta)\in C^{(n)}_\mathbb{P}$. Let us now check that  $\mathbb{P}_{j(\delta)}=\mathbb{P}\cap V_{j(\delta)}$ and that $j(\delta)$ is $\mathbb{P}$-reflecting.


For each $\alpha<j(\delta)$ the formula
$$\exists \beta\,\exists X\,(\beta>\alpha\,\wedge\, X=V_\beta\,\wedge\, \mathbb{P}_\alpha\s V_\beta)$$
is $\Sigma_{m+1}$, with parameter $\alpha$.  Since $j(\delta)\in C^{(n+m)}$, for each $\alpha<j(\delta)$ there is $\alpha<\beta<j(\delta)$ witnessing the above.  Combining this with the $\Delta_{m+1}$-definability of $\mathbb{P}$ and with $1\leq m\leq n$ it follows that 
$$\mathbb{P}_{j(\delta)}=\mathbb{P}\cap V_{j(\delta)}=\mathbb{P}^{V_{j(\delta)}}.$$ 

Let us now prove that $j(\delta)$ is $\mathbb{P}$-reflecting by showing that all the assumptions of Proposition \ref{firstprop} are met. First, $j(\delta)$ is inaccessible (actually, Mahlo) and by the above displayed expression $\mathbb{P}$ is a forcing iteration such that $\mathbb{P}_{j(\delta)}\s V_{j(\delta)}$. Second, thanks to the Mahloness of $j(\delta)$, the iteration $\mathbb{P}_{j(\delta)}$ is $j(\delta)$-cc and so it is not hard to  show that $\mathbb{P}_{j(\delta)}$ preserves the inaccessibility of $j(\delta)$.  Finally, we claim that $j(\delta)\in C(\mathbb{P})$, and so that $\forces_{\mathbb{P}_{j(\delta)}}\text{$``\mathbb{P}_{[j(\delta),\ord)}$ is $j(\delta)$-distributive''}.$

Indeed,  arguing as in the proof of Theorem \ref{preservationCn} and using that the cardinal $j(\delta)$ is $\Sigma_{m+n}$-correct, hence $\Sigma_{m+2}$-correct, one infers that $j(\delta)$ is an accumulation point of $C(\mathbb{P})$, hence $j(\delta)\in C(\mathbb{P})$.\footnote{Note that here we have used that $n\geq 2$.} 
\end{proof}
This shows that $\delta$ is $\mathbb{P}$-\Cn-extendible with $\mathbb{P}$-$\Sigma_n$-reflecting target.
\end{proof}
For simplicity let us assume that $n\geq 2$, as the argument for $n=1$ is the same. The above claim alongside the previous comments imply the existence of a proper class of $\mathbb{P}$-\Cn-extendible cardinals with $\mathbb{P}$-$\Sigma_n$-reflecting target, hence Theorem \ref{NonHomoVopenka} yields the existence of a proper class of \Cn-extendible cardinals in $V^\mathbb{P}$ and so $V^\mathbb{P}\models \mathrm{VP}(\mathbf{\Pi_{n+1}})$.
\end{proof}

\subsection{On $\rm{V=HOD}$ and the Ground Axiom}\label{SectionVHodExten}


The first forcing iteration producing a generic extension where $\rm{V}=\rm{HOD}$ holds was defined by A. McAloon \cite{McAloon}. The idea is to code the universe into the power-set function pattern so that all sets become  definable using ordinals as parameters. For more sophisticated  codings  see \cite{Broo2}.  

For the purposes of the current section we may also assume that the GCH holds, for otherwise we can force it while preserving $C^{(n)}$-extendible cardinals (cf. Theorem \ref{preserveCnext}).


 Let $\mathbb{P}$ be the class forcing notion from \cite[Theorem 3.3]{Rie}. \label{HodExtendibles} It is easy to see that $\mathbb{P}$ is $\Delta_2$-definable and that  the class of $\mathbb{P}$-reflecting cardinals $\kappa$ such that $ \mathbb{P}\cap V_\kappa=\mathbb{P}_\kappa$   contains all Mahlo cardinals, if there are any. Moreover, if there is a \Cn-extendible cardinal, there are class many Mahlo cardinals, hence the above class is proper, and every \Cn-extendible cardinal is also $\mathbb{P}$-\Cn-extendible with  $\mathbb{P}$-$\Sigma_n$-reflecting target (see Proposition \ref{PredicateandCnext}). Thus $\mathbb{P}$ fulfils the assumptions of Theorem \ref{NonHomoVopenka}.  
As the GCH holds, $\mathbb{P}$ yields a cardinal-preserving generic extension in which the  \emph{Continuum Coding Axiom} (CCA) holds, hence where  $\rm{V}=\rm{HOD}$ holds \cite[Theorem 3.3]{Rie}. Altogether, we obtain the following, which extends \cite[Theorem 3.9]{Rie}. 


\begin{cor}\label{cnextcca}
Forcing with $\mathbb{P}$ produces a generic extension of $\rm{CCA}+\neg\mathrm{GCH}$ and preserves $C^{(n)}$-extendible cardinals, for $n\geq 1$. In particular, $C^{(n)}$-extendible cardinals are consistent with $\mathrm{V=HOD}$.
\end{cor}


Likewise, the following extends \cite[Corollary 4]{HRW}:


\begin{cor}\label{cnextega}
There is a class iteration forcing  $``\mathrm{V\neq HOD}+\mathrm{GA}$'' and preserving \Cn-extendible cardinals, for $n\geq 1$. 
\end{cor}
\begin{proof}
Let $\mathbb{P}$ be the forcing iteration of Corollary \ref{cnextcca} and let $\dot{\mathbb{Q}}$ be a $\mathbb{P}$-name for the   iteration with Easton support that forces with $\mathrm{Add}(\kappa,1)$ at each regular cardinal $\kappa$ such that $2^{<\kappa}=\kappa$. Set $\mathbb{R}:=\mathbb{P}\ast\dot{\mathbb{Q}}$. By the argument in \cite[Theorem 3]{HRW}, $V^{\mathbb{R}}\models \text{$``\mathrm{V\neq HOD}+\mathrm{GA}$''}$. Since $\mathbb{R}$ is $\Delta_2$-definable and the class of $\mathbb{R}$-reflecting cardinals $\kappa$ such that $ \mathbb{R}\cap V_\kappa=\mathbb{R}_\kappa$ contains all Mahlo cardinals,  every $C^{(n)}$-extendible cardinal is $\mathbb{P}$-$C^{(n)}$-extendible with $\mathbb{P}$-$\Sigma_n$-reflecting target.  As  $\mathbb{R}$ satisfies the hypotheses  of  \ref{NonHomoVopenka},   the result follows.  
\end{proof}

\bibliographystyle{alpha} 
\bibliography{biblio-3}
\end{document}